\documentclass[12pt]{article}
\usepackage{amsmath,amsfonts,amssymb,amsthm, graphicx, epsfig, float}
\usepackage{listings,paralist}

\newcommand{\I}{{\bf 1}}
\newtheorem{proposition}{Proposition}[section]
\newtheorem{theorem}[proposition]{Theorem}

\newtheorem{lemma}[proposition]{Lemma}
\newtheorem{remark}[proposition]{Remark}

\numberwithin{equation}{section}
\newcommand{\nc}{\newcommand}
\nc{\R}{{\mathbb R}}
\nc{\bS}{{\mathbb S}^{d-1}}
\nc{\N}{{\mathbb N}}
\nc{\Z}{{\mathbb Z}}
\nc{\BP}{\mathbb{P}}
\nc{\BE}{\mathbb{E}}
\nc{\BQ}{\mathbb{Q}}
\nc{\bN}{{\mathbf N}}
\nc{\cK}{{\mathcal K}^d}
\nc{\cN}{{\mathcal N}}
\DeclareMathOperator{\BV}{{\mathbb Var}}
\DeclareMathOperator{\card}{card}

\begin{document}

\author{Sven Ebert\footnote{
Institut f\"ur Stochastik,
 Karlsruher Institut f\"ur Technologie,
76128 Karlsruhe, Germany.
Email: sven.ebert@kit.edu}
\ and
G\"unter Last\footnote{
Institut f\"ur Stochastik,
 Karlsruher Institut f\"ur Technologie,
76128 Karlsruhe, Germany.
Email: guenter.last@kit.edu}
}

\title{On a class of growth-maximal hard-core processes}
\date{\today}
\maketitle

\begin{abstract}
\noindent
Generalizing the well-known lilypond model (\cite{HM96,DSS99,DL05})
we introduce a growth-maximal hard-core model based on a space-time point process
$\Psi$ of convex particles. Using a purely deterministic algorithm we prove
under fairly general assumptions that the model exists and
is uniquely determined by $\Psi$. Under an additional stationarity assumption
we show that the model does not percolate. Our model generalizes
the lilypond model considerably even if all grains are born at the same time.
In that case and under a Poisson assumption we prove a central limit theorem
in a large volume scenario.
\end{abstract}



\section{Introduction}\label{secintro}

We consider a point process $\Psi=\{(X_n,T_n,Z_n):n\ge 1\}$
on $\R^d\times\R_+\times\cK$, where $\R_+:=[0,\infty)$
and $\cK$ denotes the space of all convex bodies containing
the origin $0\in\R^d$ in its interior.
We interpret $X_n$ as the position of a grain (particle)
with shape $Z_n$ that is born at time $T_n$.
Without interaction the $n$-th grain
starts growing at time $T_n$ to form a grain $X_n+(t-T_n)Z_n:=\{X_n+(t-T_n)x: x\in Z_n\}$
at time $t\ge T_n$. It ceases its growth as soon as it encounters
any other particle. This interaction can occur in two ways.
First,  $X_n$ can be covered by another grain by time $T_n$, so that
no growth can start at all. Second, the growing grain may touch
another grain (which may itself be either growing or have
ceased growing at an earlier time).
As a result of this growth process a growth
time $R_n\in\R_+$ is attached to the $n$-th particle. The full-grown
grain is given by $Z_n^*:=X_n+R_nZ_n$ if $R_n>0$.
We use this notation also in case $R_n=0$. (Then $Z^*_n=\{X_n\}$
has empty interior.) The growth protocol (informally) described
above entails the following three properties.
\begin{enumerate}
\item[(i)]
The interiors of $Z_m^*$ and $Z_n^*$ do not intersect whenever
$m\ne n$.
\item[(ii)] If $R_m=0$ then there is some $n\ne m$ such
that $R_n>0$, $X_m\in Z^*_n$ and
$T_m\ge T_n+\inf\{r\ge 0:X_m\in X_n+rY_n\}$.
\item[(iii)] If $R_m>0$ then there is some $n\ne m$ such
that $R_n>0$, $T_m+R_m\ge T_n+R_n$ and $Z^*_m\cap Z^*_n\ne\emptyset$.
\end{enumerate}
Property (i) says that $\{Z_n^*:R_n>0\}$ is a {\em hard-core}
system, while (ii) means that the growth of the $n$-th particle
is inhibited by some other grain. Property (iii) says
that  any grain with positive growth time is stopped
by some other grain that reached its final size
at the same time or earlier. We summarize
properties (ii) and (iii) by calling
$(X_n,T_n,Z_n)$ (or $Z_n^*$) an {\em earlier neighbour}
of $(X_m,T_m,Z_m)$ (resp.\ $Z_m^*$). Any grain must have at least one
earlier neighbour.
If $\{R_n:n\ge 1\}$ is a family of $\R_+$-valued random variables
depending measurably on $\Psi$ and such that
(i)-(iii) hold almost surely, then we call
$\Psi^*=\{(X_n,T_n,Z_n,R_n):n\ge 1\}$
a {\em growth-maximal hard-core model} (based on $\Psi$).
Under certain assumptions on $\Psi$ we will prove
in Section \ref{secexist} that the model exists
and is uniquely determined by $\Psi$.

\begin{figure}[h]
\begin{center}
\begin{minipage}[c]{5in}
\parbox[c]{2.5in}{
\includegraphics[width=2.4in]{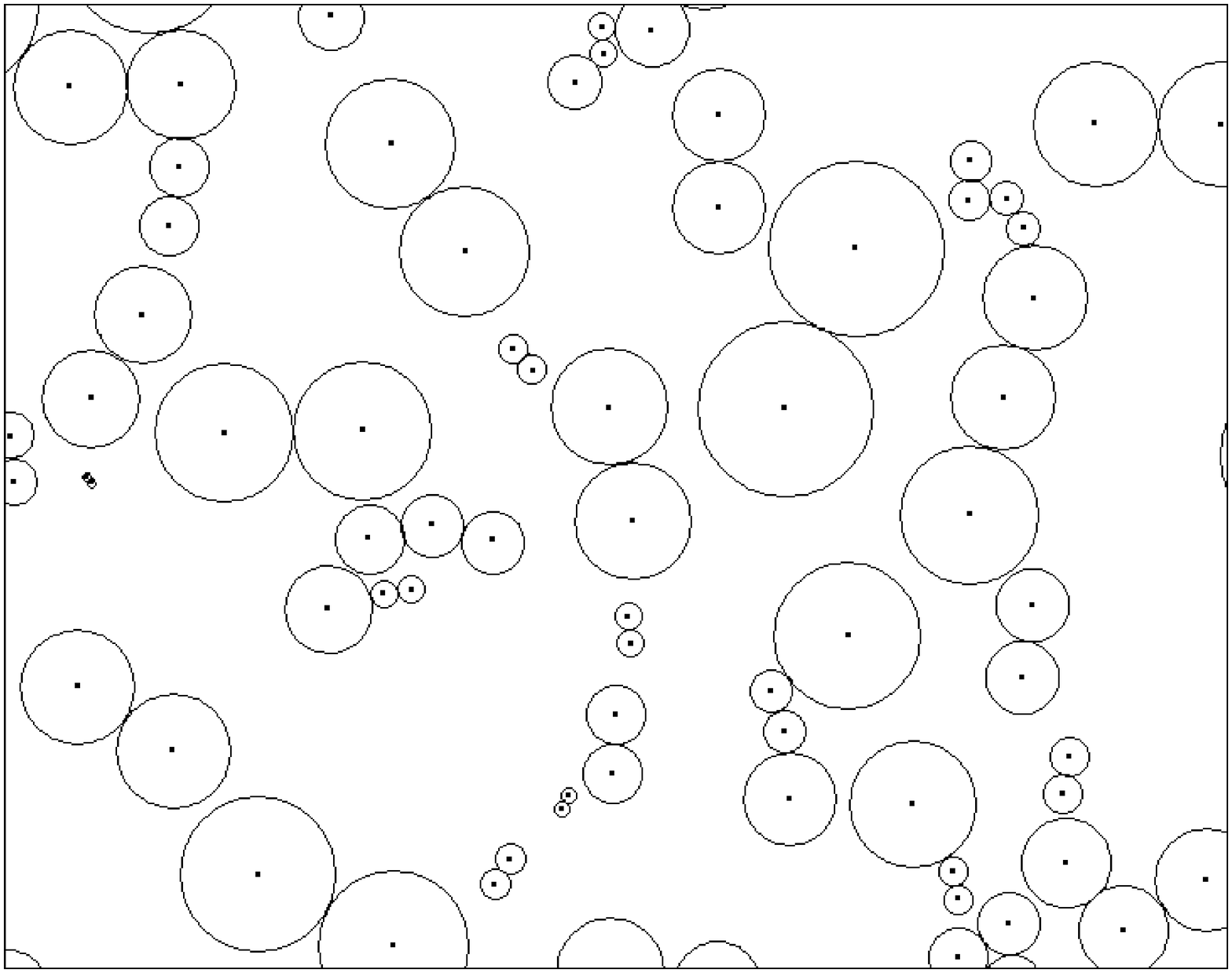}
} \hfill
\parbox[c]{2.5in}{
\includegraphics[width=2.4in]{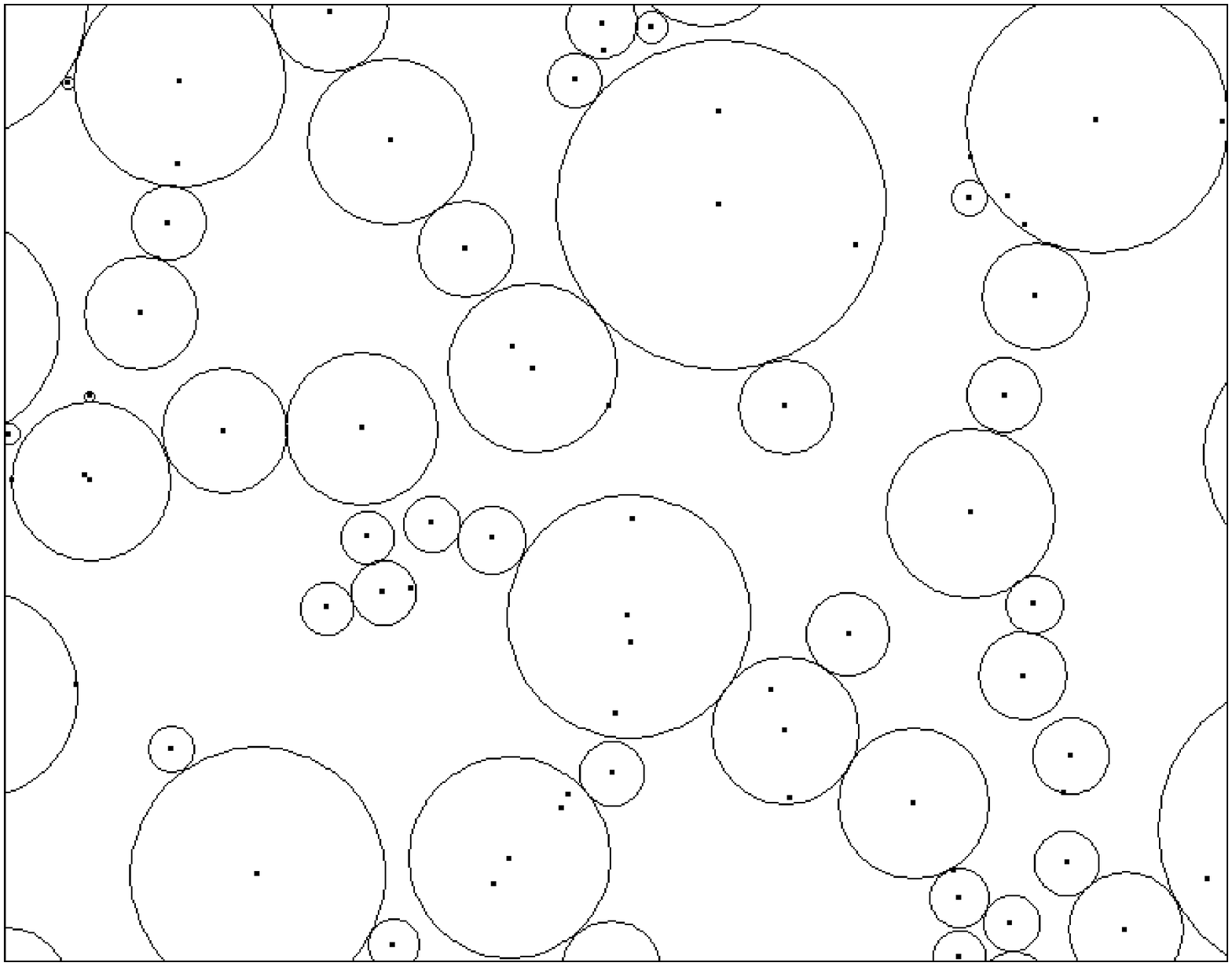}
}
\end{minipage}
\end{center}
\caption{\label{compdynclass}
Two models based on the same configuration of positions but with different birth time scenarios.
}
\end{figure}

The {\em lilypond model} \cite{HM96,DSS99,DL05} is a well-known
special case of a growth-maximal hard-core model. In this
case $(T_n,Z_n)=(0,B^d)$ for all $n\ge 1$, where $B^d$
is the (closed) unit ball in $\R^d$. In the physics literature this was
called {\em touch-and-stop model}, see \cite{AnBrKr94}.
In \cite{HL06} the model was generalized to the case $(T_n,Z_n)=(0,B)$ for a
general symmetric star-body $B\subset\R^d$.
Neither the case of random (possibly different) shapes $Z_n$ nor the case
of a space-time driving process $\Psi$ has been studied before.

\begin{figure}[h]
\begin{center}
\begin{minipage}[c]{4in}
\frame{
\includegraphics[width=4in]{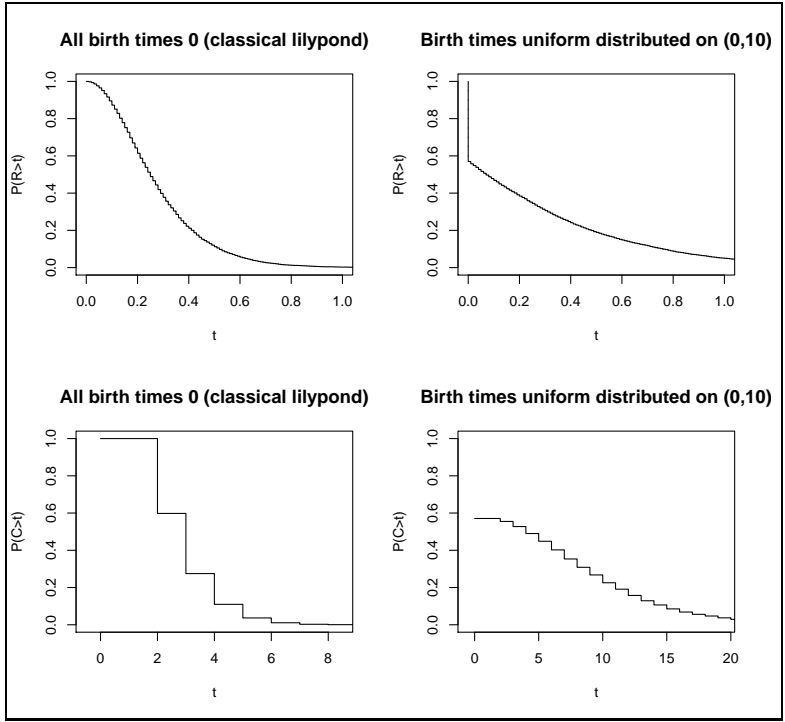}
}
\end{minipage}
\end{center}
\caption{\label{lengthclustercomp}
Comparison of length and cluster distributions in different birth time scenarios.
}
\end{figure}

Figure \ref{compdynclass} gives an impression of how the
introduction of starting times alters the structure of the model.
The set of grains is the same in both pictures. But on the left-hand side all
birth times are zero while on the right-hand side the birth times are independent
and uniformly distributed on the interval $[0,10]$.
We see that in the presence of birth times some germs might not be born at all.
This in turn implies that the germs which are born have more space to grow.
Figure \ref{lengthclustercomp} compares the length and cluster size distributions
of the preceding two birth time scenarios. The distributions on the
right-hand side have considerably
higher probabilities at large values.

The paper is organized as follows. In Section \ref{secexist}
we use an exteneded version of an algorithm in \cite{DL05}
to prove under fairly general assumptions that the model exists
and is uniquely determined by $\Psi$ in a translation invariant way.
In Section \ref{secperc} we generalize the results in \cite{HM96,DL05}
and prove under appropriate stationarity and
moment assumptions that there is no percolation in the model.
In Section \ref{secstabclt} we consider an independently
marked stationary Poisson process with constant birth times.
Using stabilization arguments as in \cite{LaPe12},
we prove that the growth times $R_n$, $n\in\N$, satisfy a central limit theorem.

\section{Existence and uniqueness}\label{secexist}

Consider $\Psi=\{(X_n,T_n,Z_n):n\ge 1\}$, as introduced above,
defined on an abstract probability space $(\Omega,\mathcal{F},\BP)$.
Assume that $\Phi:=\{X_n:n\ge 1\}$ is {\em locally finite},
that is $\Phi(B):=\card\{m\ge 1:X_m\in B\}$ is $\BP$-a.s.\ finite
for compact sets $B\subset\R^d$.
For any $k\ge 1$, let $\nu^{(k)}$ denote the
{\em $k$-th factorial moment measure} of $\Psi$,
the measure on $(\R^d\times\R_+\times\cK)^k$ defined
by
\begin{align}\label{deffacmommeas}
\nu^{(k)}(C):=\BE\sum\nolimits^*\I\{(x_1,t_1,K_1,\dots,x_k,t_k,K_k)\in C\cap \Phi^k\}
\end{align}
for measurable sets $C\subset (\R^d\times\R_+\times\cK)^k$,
where the * indicates that the summation is over
all $k$-tuples $((x_1,t_1,K_1),\dots,(x_k,t_k,K_k))$ with
pairwise different entries. We shall assume that
\begin{align}\label{assfac}
\nu^{(k)}\le a^k(\lambda^d\otimes\BQ)^k,\quad k\ge 1,
\end{align}
where $a> 0$, $\lambda^d$ denotes the Lebesgue measure
on $\R^d$ and $\BQ$ a probability measure on
$\R_+\times\cK$. We shall further suppose that the
measure $\BQ$ satisfies the integrability condition
\begin{align}\label{assmoment}
\int (r^{2d}+1)(\rho(K)^{2d}+1)\BQ(d(r,K))<\infty,
\end{align}
where $\rho(K)$ is the radius of the smallest ball circumscribing $K\in\cK$.

In this section we will prove the following
existence and uniqueness result.

\begin{theorem}\label{texist} Assume that \eqref{assfac} and
\eqref{assmoment} hold. Then there exists a unique
growth-maximal hard-core model $\Psi^*=\{(X_n,T_n,Z_n,R_n):n\ge 1\}$
based on $\Psi$.
\end{theorem}
Under the given assumptions the theorem says that there is
a family $\{R_n:n\ge 1\}$ of $\Psi$-measurable random
variables such that $\Psi^*$ satisfies conditions (i)-(iii) given
in the introduction almost surely.
Moreover, if $\{R'_n:n\ge 1\}$ is another such family,
then $R_n=R'_n$ $\BP$-a.s.\ for any $n\ge 1$.

Before proving Theorem \ref{texist} we provide a short
discussion of assumption \eqref{assfac}. For $k\ge 1$ let
$\mu^{(k)}$ denote the $k$-th factorial moment measure of
$\Phi:=\{X_n:n\ge 1\}$. If $\Psi$ is an
{\em independent marking} of $\Phi$ with mark distribution
$\BQ$, then
$$
\nu^{(k)}(d(x_1,t_1,K_1,\dots,x_k,t_k,K_k))=
\mu^{(k)}(d(x_1,\dots,x_k))\BQ^k(d(t_1,K_1,\dots,t_k,K_k))
$$
and \eqref{assfac} is implied by $\mu^{(k)}\le a^k (\lambda^d)^k$,
$k\ge 1$. The most fundamental example of a point process $\Phi$
with this property is a Poisson process with a bounded intensity
function (e.g.\ a stationary Poisson process). Other examples are
provided by certain Cox, Poisson cluster and Gibbs processes, see
\cite{DL05} for more details.

Our proof of Theorem \ref{texist} is based on a purely
deterministic result inspired by
Proposition 2.8 in \cite{DL05}. We start by introducing
some notation. Define $\bN$ as the set of
all countable sets $\psi\subset \R^d\times\R_+\times\cK$
such that
$\card\{v\in\psi :v\in B\times\R_+\times\cK\}<\infty$ for all bounded sets $B\subset\R^d$
and for all $(x,s,K),(y,t,L)\in\psi$ we have
$(x,s,K)\ne (y,t,L)$ only if $x\ne y$.
Any $\psi\in\bN$ is identified with the (counting) measure
$\card\{v\in \psi:v \in \cdot\}$ on $\R^d\times\R_+\times\cK$.
A point $(x,s,K)\in\psi$ is interpreted as a {\em grain} with
{\em shape} $K$ located at $x$ and born at time $s$.
As usual we equip $\bN$ with the smallest
$\sigma$-field $\mathcal{N}$ making the mappings $\psi\mapsto\psi(C)$ measurable
for all measurable $C\subset \R^d\times\R_+\times\cK$.
The point process $\Psi$ is then formally defined as a random element in $\bN$.

Consider $\psi=\{(x_n,t_n,K_n):n\ge 1\}\in\bN$. We call a subset
$\psi^*=\{(x_n,t_n,K_n,r_n):n\ge 1\}$ of
$\R^d\times\R_+\times\cK\times\R_+$ with projection
$\psi$ on $\R^d\times\R_+\times\cK$ a {\em growth-maximal hard-core
model based on $\psi$} if the interiors of $x_m+r_mK_m$ and
$x_n+r_nK_n$ do not intersect whenever $m\ne n$ and any point in
$\psi^*$ has at least one {\em earlier neighbour}: If
$r_m=0$ there is some $n\ne m$ such that $r_n>0$, $x_m\in
x_n+r_nK_n$ and $t_m\ge t_n+\inf\{r\ge 0:x_m\in x_n+rK_n\}$. If
$r_m>0$ there is some $n\ne m$ such that $r_n>0$, $t_m+r_m\ge
t_n+r_n$ and $(x_m+r_mK_m)\cap(x_n+r_nK_n) \ne\emptyset$. For
an earlier neighbour $(x_n,t_n,K_n)$ of $(x_m,t_m,K_m)$ we will
equivalently use the term {\em stopping neighbour} because of
the geometric interpretation of our model given in the introduction. Furthermore we will
refer to the function $R:\psi\rightarrow\R_+$ defined by
$R(x_n,t_n,K_n):=r_n$, $n\in\N$, as a {\em hard-core function} on $\psi$ if
$\psi^*$ is a hard-core model. If $\psi^*$ is moreover growth maximal we call
$R$ {\em a growth-maximal hard-core function}.
If $R$ is a hard-core function on $\psi$ we call two different points
$u=(x,s,K), v=(y,t,L)\in\psi$ {\em neighbours} w.r.t.\  $R$ if
$(x+R(u)K)\cap (y+R(v)L)\ne\emptyset$.

For convex and non-empty subsets $K,K',L,L'$ of $\R^d$
we define
\begin{align}\label{a}
a(K,K',L,L'):=\inf\{r\ge 0: (K+rK')\cap(L+rL')\ne\emptyset\}
\end{align}
and use abbreviations of the type $a(x,K,y,L):=a(\{x\},K,\{y\},L)$, $x,y\in\R^d$. For
two grains $u:=(x,s,K)$ and $v:=(y,t,L)$ we define the {\em first
contact time} $d(u,v)$ as follows. If $s\le t$ we set
\begin{align}\label{d}\notag
d(u,v)&:=\I\{a(x,K,y,0)\le t-s\}(s+a(x,K,y,0))\\
&\ \quad+\I\{a(x,K,y,0)> t-s\}(t+a(x+(t-s)K,K,y,L)).
\end{align}
In the case covered by the first summand,  the growing grain
$u$ reaches the point $y$ not later than the birth time $t$ of grain
$v$. In the second case both grains start to grow and meet at time $d(u,v)$.
If $s\ge t$ we define $d(u,v):=d(v,u)$, so that $d(\cdot,\cdot)$
becomes symmetric. In the following and later we write
$a\wedge b:=\min\{a,b\}$ and $a\vee b:=\max\{a,b\}$ for $a,b\in\R$.

\begin{lemma}\label{l2.1}
Let $R$ be a hard-core function on $\psi\in\bN$.
In addition suppose $u=(x,s,K)$ and
$v=(y,t,L)$ to be neighbours in $\psi$ w.r.t.\ $R$ and $R(u)\wedge R(v)>0$ . Then
\begin{align*}
(R(u)+s)\wedge (R(v)+t) \le d(u,v)\le (R(u)+s)\vee (R(v)+t).
\end{align*}
Here the second inequality is strict, provided that $R(u)+s\ne R(v)+t$.
\end{lemma}
{\sc Proof:} We assume w.l.o.g.\ $s\le t$ and prove the
assertion by contradiction.
First suppose that $(R(u)+s)\wedge( R(v)+t)> d(u,v)$.
If $a(x,K,y,0)\le t-s$ then
$$
R(u)+s>d(u,v)=s+a(x,K,y,0),
$$
where we have used \eqref{d} to get the identity.
Since $y\in x+a(x,K,y,0)K$, this contradicts $R(v)>0$ and
the hard-core property.
If $a(x,K,y,0) > t-s$, then $d(u,v)=t+r$, where $r:=a(x+(t-s)K,K,y,L)$.
Therefore $R(v)>r$. But the definition of $a(\cdot)$ implies
$(x+(t-s+r)K)\cap (y+rL)\ne\emptyset$ and hence $R(u)<t-s+r=d(u,v)-s$,
by the hard-core property. This contradicts our assumption $R(u)+s>d(u,v)$.

Now suppose $(R(u)+s)\vee (R(v)+t)< d(u,v)$.
If $a(x,K,y,0)\le t-s$, then
$$
R(v)+t<d(u,v)=s+a(x,K,y,0)\le t,
$$
contradicting $R(v)>0$. If $a(x,K,y,0)> t-s$ we obtain
$$
R(v)+t<d(u,v)=t+a(x+(t-s)K,K,y,L)=t+r,
$$
so that the definition of $a(\cdot)$ implies
$(x+(t-s+r)K)\cap (y+R(v)L)=\emptyset$. Since $u$ and $v$ are
neighbours we get $R(u)>t-s+r=d(u,v)-s$, contradicting
our assumption $R(u)+s<d(u,v)$.

To prove the final assertion we again argue by contradiction.
As the case $a(x,K,y,0)\le t-s$ can be treated as above
we assume that $a(x,K,y,0)> t-s$. Assume first that $d(u,v)\ge R(u)+s>R(v)+t$.  
As before we then get $R(u)>d(u,v)-s$, a contradiction.
Let us finally assume that $d(u,v)\ge R(v)+t>R(u)+s$. 
Then the definition of $a(\cdot,\cdot)$ implies that
$R(u)\ge t-s+r=d(u,v)-s$, contradicting 
$d(u,v)>R(u)+s$. \qed

\bigskip


The {\em translation} of $\psi\in\bN$ by $y\in\R^d$ is defined by
$y+\psi:=\{(y+x,s,K):(x,s,K)\in\psi\}$. A set $H\subset\bN$ is
{\em translation invariant} if $y+\psi\in H$ for all $\psi\in H$ and all $y\in \R^d$.
For a mapping $R:\bN\times \R^d\times\R_+\times\cK\rightarrow\R_+$
and $\psi\in\bN$ we define $R^\dag(\psi)\subset \R^d\times\R_+\times\cK\times\R_+$
by
\begin{align*}
R^\dag(\psi):=\{(x,t,K,R(\psi,x,s,K)):(x,s,K)\in\psi\}
\end{align*}
If $R$ is {\em translation invariant} in the sense that
\begin{align}\label{ti}
R(y+\psi,x+y,t,K)=R(\psi,x,s,K),\quad \psi\in\bN,\ x,y\in\R^d,\ (s,K)\in\R_+\times\cK,
\end{align}
then $R^\dag(y+\psi)=R^\dag(\psi)$ for all $(\psi,y)\in\bN\times\R^d$.

Theorem \ref{texist} is implied by the following Proposition and Lemma \ref{l2.3}.

\begin{proposition}\label{p2.2} There exist a measurable mapping
$R:\bN\times \R^d\times\R_+\times\cK\rightarrow\R_+$ and a
set $H\in\mathcal{N}$ such that if $\psi\in H$ then $R^\dag(\psi)$
is the unique growth-maximal hard-core model
based on $\psi$. Moreover $H$ and $R$ are translation invariant.
\end{proposition}

The next lemma shows that $H$ is reasonably big.

\begin{lemma}\label{l2.3}
Let $H$ be the set from Proposition \ref{p2.2}. Then $\BP(\Psi\in H)=1$.
\end{lemma}

{\sc Proof of Proposition \ref{p2.2}:} Our strategy for the proof is
as follows. We first construct $H$ and $R$ via a
recursive {\em algorithm}. Growth-maximality and the uniqueness
property will be proved in two separate lemmas.

Take $\psi\in\bN$. Starting with $\psi_0:=\emptyset$
we will recursively construct an increasing sequence $\psi_i\subset\psi$,
$i\in \N$, together with a sequence of hard-core functions $R_i:\psi_i\rightarrow\R_+$
such that $R_{i+1}$ extends $R_i$ for any $i\ge 0$.
We will further define $R(\psi,u):=R_i(u)$ whenever $u\in\psi_i$ for some $i$. Afterwards
this definition will be extended to all other $u\in\psi$ in a suitable way.

Let $i\in\N$. For $u\in\psi\setminus\psi_i$
we define
\begin{align*}
\nu_i(u):=\{u\}\cup\{v\in\psi\setminus(\psi_i\cup\{u\}):&\ v \mbox{
and } u \mbox{ are mutual nearest neighbours}\\
& \mbox{ in } \psi\setminus\psi_i \mbox{ with respect to }
d\}
\end{align*}

We assume from now on that $\psi$ satisfies
\begin{align}\label{c2}
d(u,v)\ne d(u',v'),\quad u,u',v,v'\in\psi,\, \{u,v\}\ne \{u',v'\}.
\end{align}
Then $\nu_i(u)$ contains at most two
points. Furthermore let
\begin{align}\label{defvarphi_i}
\varphi_i:=\{u\in\psi\setminus\psi_i: \nu_i(u)\neq\{u\}\}
\end{align}
be the set of all mutual nearest neighbours in
$\psi\setminus\psi_i$ (w.r.t. $d$). In the case $\varphi_i=\emptyset$ we
define $\psi_{i+1}=\psi_i$ i.e. $\psi_{j}=\psi_i$, $j\ge i$. Otherwise
we obtain $\psi_{i+1}$ based on $\psi_i$ in the following way.
Set $f(u,v):=\inf\{r\ge 0: (x+ r K)\cap (y + R_i(v)L)\neq\emptyset\}+s$ for
$u=(x,s,K)\in\psi\setminus\psi_i$ and $v=(y,t,L)\in\psi_i$.
For $u\in\varphi_i$ define $d_i(u):=\inf_{v\in\psi_i,R_i(v)>0} f(u,v)$ as well as
\begin{align*}
s'_i(u):=d(u,v),\ v\in\nu_n(u)\setminus\{u\}\quad \mbox{and}\quad
s''_i(u):=\min_{v\in\nu_i(u)} d_i(v).
\end{align*}
Now we add points of $\varphi_i$ to
$\psi_i$ by the following procedure to obtain $\psi_{i+1}$.

Consider $u=(x,s,K)\in\varphi_i$ so that there exists
$v=(y,t,L)\in\nu_i(u)\setminus\{u\}$. Assume first that $s'_i(u)\le s''_i(u)$.
If $s\le t$ and $a(x,K,y,\{0\})\le t-s$ we set
$R_{i+1}(v):=0$ and add $v$ to $\psi_i$. If $s\le t$ and
$a(x,K,y,\{0\})> t-s$ we add $u$ and $v$ to $\varphi_i$ and
define $R_{i+1}(u):=s'_i(u)-s$ and $R_{i+1}(v):=s'_i(u)-t$. In case $t\le s$ we apply
the same rules as before with the roles of $u$ and $v$ interchanged.
Now we assume $s''_i(u)< s_i'(u)$. For all $w:=(z,r,M)$ with
$w\in\nu_i(u)$ and $d_i(w)=s''_i(u)$, we set
$R_{i+1}(w):=s''_i(u)-r$ and add those points to $\psi_i$. Applying
this algorithm to every $u\in\varphi_i$ yields the $(i+1)$-th
hard-core function $R_{i+1}$ on $\psi_{i+1}$.

We define
\begin{align}\label{defvarphi}
\psi_\infty:=\bigcup_{i=1}^\infty\psi_i\quad \mbox{ and }\quad
\varphi:=\psi\setminus\psi_\infty=\bigcap_{i=1}^\infty\psi\setminus\psi_i.
\end{align}
Set $R(\psi,u):=R_i(u)$ whenever $u\in\psi_i$ for some $i$.
If $\card \varphi =1$ and $\{v\in\psi_\infty:R(v)>0\}\ne\emptyset$
we define
$R(\psi,u):=\inf_{v\in\psi_\infty,R(\psi,v)>0}f(u,v)-s$ for $u=(x,s,K)\in\varphi$.

Let $H$ denote the set of all $\psi$ satisfying \eqref{c2} and for which the
above algorithm yields a set $\psi_\infty$ with $\card (\psi\setminus \psi_\infty)\le 1$ and
$\{v\in\psi_\infty:R(\psi,v)>0\}\ne\emptyset$. For all $\psi\in\bN\setminus H$,
we set $R(\psi,\cdot)\equiv0$. Then $R$ is a hard-core function.
By construction the set $H$ is translation invariant,
and $R$ is translation-invariant on $H$. Thus $R$ is translation-invariant. 
We need to show that $R(\psi,\cdot)$ is for any $\psi\in H$ the unique growth-maximal
hard-core function on $\psi$. To aid reading, we do this in two separate
Lemmas. 

\begin{lemma}\label{lgrowthm} Assume that $\psi\in H$. Then
$R:=R(\psi,\cdot)$ is a growth-maximal function.
\end{lemma}
{\sc Proof:} Let $u=(x,s,K)\in\psi$. To show that $u$ has an earlier neighbour.
We distinguish the cases $R(u)=0$ and $R(u)>0$.

1. Assume that $R(u)=0$. By the algorithm (introduced in the above first part
of the proof of  Proposition \ref{p2.2}) there is an $u\in\psi_\infty$, that is 
$u\in\psi_{i+1}\setminus\psi_i$ for some $i\ge 0$. Moreover, the algorithm implies
that there exists a $v=(y,t,L)\in\psi\setminus\{u\}$ such that
$\{u,v\}=\nu_i(u)$ and 
\begin{align}\label{++}
s'_i(v)\le s''_i(v).
\end{align} 
We show, that $v$ is an earlier neighbour of $u$. First we assume that
$v\in\psi_{j+1}\setminus\psi_j$ for some $j>i$ 
(the other case being $v\in\psi\setminus\psi_\infty$).
Thus there is a $w=(z,r,M)\in \psi\setminus(\psi_j\cup\{u,v\})$ with $\{v,w\}=\nu_j(v)$. 

Assume now that $s'_j(v)\le s''_j(v)$. Since $u,v$ are mutual nearest neighbours in
$\psi\setminus\psi_i\supset \psi\setminus\psi_i$ and \eqref{c2} holds 
(since $\psi\in H$), we have
\begin{align}\label{*}
d(v,w)> d(v,u).
\end{align}
We now claim $R(v)>0$. Assuming on the contrary that $R(v)=0$ we obtain from
$R(u)=0$ and \eqref{d} that $t\le s$, 
\begin{align}\label{+}
a(x,0,y,L)=\le s-t
\end{align}
and
\begin{align}\label{**}
t+a(x,0,y,L)=d(u,v).
\end{align}
Using the same argument for the pair $(v,w)$ we obtain
\begin{align*}
t+a(x,0,y,L)=d(u,v)< d(w,v) = r+a(y,0,z,M)
\end{align*}
and, moreover, $a(y,0,z,M)\le t-r$. It follows that
$t+a(x,0,y,L)<t$, a contradiction. Recalling $s'_j(v)\le s''_j(v)$ we derive from
$R(v)>0$ the first equality in
$$
R(v)=d(v,w)-t>d(v,u)-t=a(x,0,y,L)
$$
and therefore
\begin{align}\label{***}
x\in y+R(v)L.
\end{align} 
Together with \eqref{+} this implies that $v$ is indeed an earlier
neighbour of $u$. 

Assume now that $s'_j(v)>s''_j(v)$. By the algorithm there exists
a neighbour $w':=(z',r',M')\in\psi_j$ of $v$ with $R(w')>0$. Further there
are $k<j$ and $w''=(z'',r'',M'')$, such
that $w'$ and $w''$ are mutual nearest neighbours in
$\psi\setminus\psi_k$. Hence the algorithm implies the first identity in
\begin{align*}
R(w')+r'=s'_k(w')\wedge s''_k(w')\le s'_k(w')=d(w',w'')\le d(w',v),
\end{align*}
where the final inequality comes from the fact that $w',w''$ are mutual nearest
neighbours in $\psi\setminus\psi_k\supset\psi\setminus\psi_j$. 
Since $w'$ and $v$ are neighbours, Lemma \ref{l2.1} yields 
\begin{align}\label{****}
R(v)+t\ge d(v,w').
\end{align} 
Again we have to distinguish two subcases.
If $k\ge i$ we get from the fact that $u$ and $v$ are mutual nearest neighbours in
$\psi\setminus\psi_i\supset\psi\setminus\psi_k$ 
that $d(w',v)\ge d(u,v)$.  Together with \eqref{****} and \eqref{**} this yields
$$
R(v)\ge d(v,w')-t\ge d(u,v)-t=a(x,0,y,L).
$$ 
Hence \eqref{***} follows. Assume now that $k<i$.  To show that
$$
R(v)=a(y,L,z'+R(w')M',0)\ge a(x,0,y,L),
$$
we assume on the contrary that the inequality fails.
But then the algorithm implies $s'_i(v)>s_i''(v)$, contradicting \eqref{++}.

Now we assume that $v\in\psi\setminus \psi_\infty$. To
show that $v$ is an earlier neighbour of $u$ we can proceed analogously to
the case $s'_j(v)> s''_j(v)$ from above.

2. Assume that $R(u)>0$. We have to consider two subcases. If $s'_i(u)\le
s''_i(u)$ there exists $v=(y,t,L)\in\psi$ such that $\{u,v\}=\nu_i(u)$, 
$R(u)=d(u,v)-s>0$ and $R(v)=d(u,v)-t>0$. Therefore
$R(u)+s=R(v)+t$ and $v$ is an earlier neighbour of $u$. If
$s'_i(u)> s''_i(u)$ then $i\ge 1$ and there exists 
a neighbour $v=(y,t,L)\in\psi_i$ of $u$ with 
$R(v)>0$. We show $R(u)+s\ge R(v)+t$ by assuming on the contrary that
$R(u)+s< R(v)+t$. There exist a $k\in\{0,\ldots,i-1\}$ and a $w$ with
$v\in\psi_{k+1}\setminus\psi_{k}$ and $\{v,w\}=\nu_i(v)$. The 
algorithm implies that $R(v)+s\le d(v,w)$. 
Lemma \ref{l2.1} and
our assumption $R(u)+s< R(v)+t$ yield $R(u)+s > d(u,v)$. This gives
$d(u,v)<d(v,w)$ which contradicts the mutual nearest neighbour
property of $v$ and $w$ in $\psi\setminus\psi_{k-1}$.

In the remaining case $u\in\psi\setminus\psi_\infty$
we can proceed as in the second case of the last
paragraph. Thus every $u\in\psi$ has at least one earlier neighbour. \qed 

\begin{lemma}\label{lgrowthuniqe} Assume that $\psi\in H$. Then
$R:=R(\psi,\cdot)$ is the unique growth-maximal function on $\psi$.
\end{lemma}
{\sc Proof:} Let $R'$ be some
growth-maximal hard-core function on $\psi$. First we prove by 
by induction over the steps of the algorithm that $R(u)=R'(u)$ 
whenever $u\in\psi_\infty$ satisfies $R(u)>0$. 

Let $u=(x,s,K) \in\psi_1$ satisfy $R(u)>0$.  The then yields algorithm that
\begin{align}\label{11} 
d(u,v)=R(u)+s=R(v)+t,
\end{align}
where $v=(y,t,L)$ is the mutually nearest neighbour of $u$ in $\psi$.
Assume $R'(u)=0$. Then there exists $w=(z,r,M)\in\psi\setminus\{u\}$ with
$R'(w)\ge a(z,M,x,0)$ and $s-r>a(z,M,x,0)$. Hence 
\begin{align}\label{12}
d(u,w)=a(z,M,x,0)+r<s.
\end{align}
If $s\le t$, we get $d(u,v)=t+a(x+(t-s)K,K,y,L)> s$. If $t<s$, we also get
$d(u,v)=s+a(x,K,y+(s-t)L,L)> s$. 
In view of \eqref{12} we obtain that $w\ne v$ and $d(u,w)<d(u,v)$.
But $w\neq v$ leads to $d(u,w)<d(u,v)$. 
As this contradicts the mutual nearest neighbour
property of $u$ and $v$, we must have that $R'(u)>0$.

Assume that $R(u)\ne R'(u)$. By Lemma \ref{l2.1} and \eqref{11}
we can assume w.l.o.g.\ that 
\begin{align}\label{13}
R'(u)+s<d(u,v).
\end{align}
Let $w\in\psi$ be a stopping neighbour of $u$ w.r.t.\ $R'$. In particular, $R'(w)>0$. 
By Lemma \ref{l2.1} we conclude $R'(u)+s\ge d(u,w)$.
Hence \eqref{13} implies that  $d(u,w)<d(u,v)$, a contradiction.

Now we take $u=(x,s,K)\in\psi_{i+1}\setminus\psi_i$  for some $i\ge 1$ and assume that 
$R(u)>0$. Again we first show that $R'(u)>0$. Assume $R'(u)=0$. Then there exists
$w\in\psi\setminus\{u\}$ with $R'(w)\ge a(z,M,x,0)>0$ and $s-r\ge a(z,M,x,0)$. Assume
$w\in\psi\setminus\psi_i$. Then we can argue as in the case $u\in\psi_1$
to conclude that $R'(u)>0$. Hence we can assume that $w\in\psi_i$. 
from the induction hypothesis we get $R'(w)=R(w)$. Therefore,
\begin{align}\label{14}
d(w,u)=a(z,M,x,0)+r<R(w)+r. 
\end{align}
Due to the algorithm $0<R(w)+r\le d(w,w')$, where $w'$ is a mutually nearest
neighbour (w.r.t.\ $R$) of $w$ in $\psi\setminus\psi_j$ for some $j\le i-1$.
Therefore $d(w,u)\le d(w,w')$ contradicting the choice of 
$w'\in \psi\setminus\psi_j\supset \psi\setminus\psi_i$ and the fact that $u\in\psi\setminus\psi_i$. 
Thus $R'(u)>0$.  

Next we show that $R(u)=R'(u)$. Assume $\nu_i(u)=\{u,v\}$ and 
$$
R(u)=s'_i(u)-s=d(u,v)-s. 
$$
Assume moreover that $R'(u)<s'_i(u)-s$. Then we can argue as at \eqref{13}
to conclude that  $u$ cannot have an earlier neighbour in $\psi\setminus\psi_i$ 
w.r.t.\ $R'$. Hence there exists an earlier neighbour $w\in\psi_i$ of $u$ w.r.t.\ $R'$,
that is  $R'(w)+r\le R'(u)+s$ and $R'(w)>0$. 
By induction hypothesis $R(w)=R'(w)$. Therefore the algorithm implies
\begin{align*}
f(u,w)=\inf\{r\ge 0:(x+rK)\cap (z+R(w)M)\ne\emptyset\}+s=R'(u)+s<s'_i(u).
\end{align*}
This 
leads to 
\begin{align*}
f(u,w)<s'_i(u)\le s''_i(u)=\inf_{v\in\psi_i,R(v)>0}f(u,v)\le f(u,w)<s'_i(u),
\end{align*}
a contradiction. Now assume $R'(u)>s'_i(u)-s$. The hard-core property of $R'$ yields
that $R'(v)<s_i(v)-t$. But we have just proved that this is not possible. 
Hence only the desired case $R'(u)=s'_i(u)-s=R(u)$ remains. 

Now assume $R(u)+s=s''_i(u)<s'_i(u)$. Then
there exists $v\in\psi_i$ with $R(v)>0$ and $s''_i(u)=f(u,v)+s$. By induction hypothesis $R(v)=R'(v)$.
Hence the hard-core property of $R'$ that
$$
s''_i(u)+s=f(u,v)\ge R'(u). 
$$
To show that $R(u)=R'(u)$ we assume on the contrary that the above inequality is strict. 
Since, by induction hypothesis, $R(w)=R'(w)$ for all $w\in\psi_i$ with
$R(w)>0$, we obtain from $s''_i(u)+s=f(u,v)$ and the definition of $f(u,v)$
that $u$ has no neighbour in $\psi_i$ w.r.t.\ $R'$.
Because of $s''_i(u)<s'_i(u)$ we get that 
$0<R'(u)+s<f(u,v)< d(u,v')$ where $v'$ is the mutual nearest neighbour of $u$ in $\psi_{i+1}\setminus\psi_i$. By Lemma \ref{l2.1} we deduce that
$$
R'(w)+r\le d(u,w)>d(u,v')>R'(u)+s,
$$
for every neighbour $w=(z,r,M)\in\psi\setminus\psi_i$ of $u$.
Therefore $U$ cannot have an earlier neighbour in $\psi\setminus\psi_i$
w.r.t.\ $R'$.
This is a contradiction.

Next we show $R(u)=R'(u)$ for $u\in\psi\setminus\psi_\infty$. 
Due to the algorithm we can assume $R(u)>0$. Because of $R(v)=R'(v)$ 
for all $v\in\psi_\infty$ with $R(v)>0$ and
the hard-core property of $R'$ we obtain the inequality
$R'(u)\le \inf_{v\in\psi_\infty, R(v)>0}f(u,v)-s=R(u)$ . 
If $R'(u)<R(u)$ the point $u$
would not have an earlier neighbour by the definition of $f(u,v)$.

Finally we show that $R'(u)=0$ holds for all $u\in\psi$ with $R(u)=0$. 
For $u\in\psi$ with $R(u)=0$ let $v\in\psi$ be an earlier neighbour of $u$ in $\psi$ 
w.r.t.\ $R$. Therefore it 
follows $R(v)>0$ and $x\in R(v)L$. We have already shown above $R'(v)=R(v)$.
Thus $x\in R'(v)L$. Since $R'$ is a hard-core function this yields $R'(u)=0$.\qed

\bigskip

For a Borel set $A\subset\R^d$  $V_d(A)$ the {\em volume} of $A$.
Further let $\kappa_d:=V_d(B^d)$ denote the volume of the unit ball.
The proof of Lemma \ref{l2.3} uses the following lemma.

\begin{lemma}\label{l2.5}
Let $r\ge 0$, $0\le s\le t$ and $K_1,K_2\in\cK$ with 
$K_1\subset x_1+r_1B^d$, $K_2\subset x_2+r_2B^d$ for some $r_1,r_2\ge 0$
and $x_1,x_2\in\R^d$.
Then 
\begin{align*}
V_d(rK_1+tK_2)-V_d(rK_1+sK_2)\le \max_{i=1,\ldots,d}(t^{i}-s^{i}) 
\sum_{j=0}^{d-1} \kappa_d\binom{d}{j} r^j r_1^j r_2^{d-j}.
\end{align*}
\end{lemma}
{\sc Proof:} Using mixed volumes (see e.g.\ \cite[(5.1.26)]{Schneider93}),
we obtain
\begin{align*}
V_d(rK_1+tK_2)=\sum_{j=0}^{d} \binom{d}{j} r^jt^{d-j}V(K_1[j], K_2[d-j]),
\end{align*}
where we have also
used the homogeneity part of \cite[(5.1.24)]{Schneider93}.
Therefore
\begin{align*}
V_d(rK_1+tK_2)-V_d(rK_1+sK_2)=\sum_{j=0}^{d-1} \binom{d}{j} r^j(t^{d-j}-s^{d-j})
V(K_1[j], K_2[d-j]).
\end{align*}
Since $K_i\subset x_i+r_iB^d$ we can now use the monotonicity property
\cite[(5.1.23)]{Schneider93} and the translation
invariance of mixed volumes to obtain that
\begin{align*}
V_d(rK_1+tK_2)-V_d(rK_1+sK_2)\le \sum_{j=0}^{d-1} \binom{d}{j} 
r^j(t^{d-j}-s^{d-j})r_1^jr_2^{d-j}V_d(B^d),
\end{align*}
where we have used that $V(B^d,\ldots,B^d)=V_d(B^d)$. This implies the
assertion.
\qed

\bigskip

%

\bigskip
{\sc Proof of Lemma \ref{l2.3}:} We construct a measurable set $G\subset H$ for which we can show that $\BP(\Psi\in G)=1$. Consider an increasing sequence $b_i>0$, $i\in\N_0$, with $\lim_{i\to\infty}b_i=\infty$ to be specified later and define
$W_j:=[-j,j]^d$, $j\ge 1$.
Moreover for $n\in\N$, $s<t$ and a bounded Borel set $B\subset\R^d$
the set $F_n(s,t,B)$ is
defined as follows. A set $\psi\in\bN$ belongs to $F_n(s,t,B)$, if it contains
$n+1$ different points $u_0,\ldots,u_n$ such that $s<d(u_0,u_1),\ldots, d(u_{n-1},u_n)\le t$
and $u_0\in B':=B\times\R_+\times \cK$. Define
\begin{align}\label{defF1}
F(s,t,B):=\bigcap_{n\ge 1} F_n(s,t,B).
\end{align}
and $G$ as follows.
A set $\psi\in\bN$ belongs to $G$ if it satisfies \eqref{c2} and
\begin{align}
\label{c1}
&\card\{v\in \psi\setminus\{u\}:d(u,v)\le r\}<\infty,\quad r>0,\,u\in\psi,\\
\label{c3} &a(x,K,y,L)\ne |s-t|,\quad (x,s,K),(y,t,L)\in\psi,\,x\ne y,
\end{align}
as well as $\psi\notin F$ where
\begin{align}\label{defF2}
F:=\bigcup_{i=0}^\infty\bigcup_{j=1}^\infty F(b_i,b_{i+1},W_j).
\end{align}
We recall the definition of the set $H$ from the proof of Proposition \ref{p2.2} as the set of all $\psi\in\bN$ which satisfy \eqref{c2} and the algorithm of Proposition \ref{p2.2} yields (based on $\psi$) a set $\psi_\infty$ with $\card (\psi\setminus \psi_\infty)\le 1$ and
$\{v\in\psi_\infty:R(\psi,v)>0\}\ne\emptyset$. We now show $G\subset H$ an take $\psi\in G$.
Define $\varphi:=\psi\setminus\psi_\infty$ and assume $\card\varphi\ge 2$.
Furthermore we take $u\in\varphi$ and define $\nu(u)$ as the
set of all mutual nearest neighbours of $u$ in $\varphi$ and set $\nu(u)=\{v\}$.
Define $A(u):=\{w\in\psi:d(u,w)\le d(u,v)\}\cup\{w\in\psi:d(v,w)\le
d(v,u)\}$. By construction $A(u)$ contains no points of
$\varphi\setminus(\nu(u)\cup\{u\})$. Because of \eqref{c1} the set
$A(u)$ is finite. So it contains no points of
$\psi_\infty$, since the corresponding growth times would
have been determined until the $i$-th step of the algorithm for some
$i\in\N_0$. Since a
pair of mutual nearest neighbours $u$ and $v$ in $\varphi$ would
also be a pair of mutual nearest neighbours in $\psi\setminus\psi_i$ at least one point of $\varphi$ would have attained a
growth time, which is a contradiction to the definition of
$\varphi$. In conclusion $\nu(u)$ has to be empty.
By induction and the symmetry of $d$ we show that in this case $\psi$
must have a descending chain, i.e.\
an infinite sequence $u_1,u_2,\dots$ of pairwise different
points of $\psi$ with $d(u_n,u_{n+1})\le d(u_{n-1},u_n)$
for any $n\ge 2$ (see also \cite{HM96} and \cite{DL05}).
Since $\varphi$ contains no pair of mutual nearest neighbours there exist
three different points $u_1,u_2,u_3$ with $d(u_1,u_2)\ge d(u_2,u_3)$ and $u_3$
being the nearest neighbour of $u_2$ and $u_2$ the nearest neighbour of $u_1$.
Inductively suppose $u_1,\ldots,u_n$, $n>2$, to be $n$ different points in $\varphi$
with $u_i$ being a nearest neighbour of $u_{i-1}$ satisfying
$d(u_i,u_{u-1})\le d(u_{i-1},u_{i-2})$, $i=3,\ldots,n$.
Let $u_{n+1}$ be a nearest neighbour of $u_n$. This implies $d(u_{n+1},u_n)\le d(u_n,u_{n-1})$
since otherwise, by the symmetry of $d$, we would get a contradiction to the
nearest neighbour property. To show $u_{n-1}\notin\{u_1,\ldots,u_n\}$ we assume
$u_{n+1}=u_j$ for some $j\in\{1,\ldots,n-1\}$. If $u_{n}$ would not be a
nearest neighbour of $u_j$ we get for some $u\in\varphi$
$$
d(u_j,u)<d(u_j,u_n)=d(u_{n+1},u_n)\le d(u_{j+1},u_j)=d(u_j,u_{j+1}).
$$
Since $u_{j+1}$ is a nearest neighbour of $u_j$ this is a contradiction and $u_{n}$
has to be a nearest neighbour of $u_j$. But this leads to $u_{n+1}$ and $u_n$ being mutual
nearest neighbours contrary to our assumption. Hence $\varphi$ would contain
a descending chain and so $\psi$ would also contain a descending chain.
This implies $\psi\in F$ (defined in \eqref{defF2}) and therefore the contradiction $\psi\notin G$.

Now assume $\psi\in G$ and
$\{v\in\psi_\infty:R(v)>0\}=\emptyset$. By our construction
algorithm this is only possible
if there exists a sequence of points $u_i=(x_i,t_i,K_i)\in\psi_\infty$,
$i\in\N_0$,  such that $t_i< d(u_i,u_{i+1})\le t_{i+1}$ and
$t_i\ge t_{i+1}$ for $i\in\N_0$. This leads to $\psi\in F$ and again to $\psi\in G^c$.

We show $\BP(\Psi\in F^c)=1$ and that the conditions 
\eqref{c2}, \eqref{c1} and \eqref{c3} have full measure, that is $\BP(\Psi\in G)=1$.
We start by calculating $\BP(\Psi\in F_n(s,t,B))$ for $n\in\N$, $s<t$ and
$B\subset\R^d$. Positive constants (depending only on the dimension) occuring 
in this calculation
are denoted by $c_i$, $i\ge 1$. Using assumption \eqref{assfac} 
and recalling the definition $B'=B\times\R_+\times \cK$ we obtain
\begin{align} \label{sum0}
\BP&(\Psi\in F_n(s,t,B))\\ \notag
&\le c_1^{n+1}\int\I\{u_0\in B'\}
\I\{s< d(u_0,u_1),\ldots,d(u_{n-1},u_n)\le t\}(\lambda^d\otimes \BQ)^{n+1}(d(u_0,\ldots,u_n)).
\end{align}
With $f(u_0,\ldots,u_{n-1}):=\I\{u_0\in B'\}\I\{s< d(u_0,u_1),\ldots,d(u_{n-2},u_{n-1})\le t\}$
we conclude
\begin{align}
  \BP&(\Psi\in F_n(s,t,B))\notag\\
\le &c_1^{n+1}\int f(u_0,\ldots,u_{n-1})\I\{s< d(u_{n-1},u_n)\le t\}
(\lambda^d\otimes \BQ)^{n+1}(d(u_0,\ldots,u_n))\notag\\
\le& c_1^{n+1}\int f(u_0,\ldots,u_{n-1})\I\{s< d(u_{n-1},u_n)\le t, t_{n-1}\le t_n\}
(\lambda^d\otimes \BQ)^{n+1}(d(u_0,\ldots,u_n))\notag\\
\label{sum}
&+c_1^{n+1}\int f(u_0,\ldots,u_{n-1})\I\{s< d(u_{n-1},u_n)\le t, t_{n-1}\ge t_n\}
(\lambda^d\otimes \BQ)^{n+1}(d(u_0,\ldots,u_n)).
\end{align}
We treat the summands separately, omitting the constant $c_1^{n+1}$
and using the definitions $\mu_n:=(\lambda^d\otimes \BQ)^{n}$,
$f_n:=f(u_0,\ldots,u_{n-1})$ and $(x_i,t_i,K_i):=u_i$, $i=0,\ldots,n$.
By \eqref{d} the first integral in \eqref{sum} is bounded by
\begin{align}\label{aaa1}
&\int f_n\I\{s<t_{n-1}+a(x_{n-1},K_{n-1},x_n,0)\le t \}\mu_{n+1}(d(u_0,\ldots,u_n))\\
\label{aaa2}
&+\int f_n\I\{s< t_n+ a(x_{n-1}+(t_n-t_{n-1})K_{n-1},K_{n-1},x_n,K_n)\le t\}
\mu_{n+1}(d(u_0,\ldots,u_n)).
\end{align}
We recall that $V_d(A)$ is the volume of $A\subset\R^d$. The integral \eqref{aaa1} equals
\begin{align}
&\int f_n\I\{t_{n-1}\le t,
(s-t_{n+1})_+<a(x_{n-1},K_{n-1},x_n,0)\le t-t_{n-1} \}\mu_{n+1}(d(u_0,\ldots,u_n))
\notag\\
&= \int f_n \I\{t_{n-1}\le t,x_n\in x_{n-1}+[(t-t_{n-1})K_{n-1}\setminus(s-t_{n-1})_+K_{n-1}]\}
\mu_{n+1}(d(u_0,\ldots,u_n))\notag\\
\label{sum11}
&\le  (t^d-s^d)\iint f_n V_d(K_{n-1})\mu_{n}(d(u_0,\ldots,u_{n-1}))\BQ(d(t_n,K_n)),
\end{align}
where we have used the inequalities $(t-c)^d-(s-c)^d\le t^d-s^d$ 
and $(t-s)^d\le t^d- s^d$ for $0\le c\le s\le t$.
We define
\begin{align}\label{defh}
h(s,t):=\max\{t^i-s^i:i=1,\ldots,d\},\quad 0\le s\le t,
\end{align}
and recall that $\rho(K)$ is the radius of the smallest all 
circumscribing $K\in\cK$. The integral \eqref{aaa2} equals
\begin{align}
\int& f_n\I\{t_{n-1}\le t,(s-t_{n})_+<a(x_{n-1}
+(t_n-t_{n-1})K_{n-1},K_{n-1},x_n,K_n)\le t-t_{n}\}\notag \\
&\qquad\quad\ \ \mu_{n+1}(d(u_0,\ldots,u_n))\notag\\
= &\int f_n \I\{t_{n-1}\le t,x_n\in x_{n-1}+(t_n-t_{n-1})K_{n-1}+(t-t_n)(K_{n-1}-K_n),\notag\\
&\qquad\quad\ \  x_n\notin x_{n-1}+(t_n-t_{n-1})K_{n-1}+(s-t_n)_+(K_{n-1}-K_n)\}
\mu_{n+1}(d(u_0,\ldots,u_n))\}\notag\\
= & \iint \I\{t_{n-1}\le t\}f_n \Big[V_d((t_n-t_{n-1})K_{n-1}+(t-t_n)(K_{n-1}-K_n))\notag\\
&\qquad\  -V_d((t_n-t_{n-1})K_{n-1}+(s-t_n)_+(K_{n-1}-K_n))\Big]
\mu_{n}(d(u_0,\ldots,u_{n-1}))\BQ(d(t_n,K_n))\notag \\ \notag
\label{sum12}\le&\iint f_n h(s,t)\sum_{i=0}^{d-1}\kappa_d\binom{d}{i}(t_n-t_{n-1})^i 
\rho(K_{n-1})^i (\rho(K_{n-1})+\rho(K_{n}))^{d-i}\\
&\qquad\quad\ \ \mu_{n}(d(u_0,\ldots,u_{n-1}))\BQ(d(t_n,K_n)),
\end{align}
where we have used Lemma \ref{l2.5} to obtain the inequality. 
Now we calculate an upper bound for the second summand in \eqref{sum}.
This can be done as in \eqref{sum11} and \eqref{sum12} with the roles of
$(t_n,K_n)$ and $(t_{n-1},K_{n-1})$ interchanged in the integrand. So we obtain
\begin{align}
\int & f(u_0,\ldots,u_{n-1})\I\{s< d(u_{n-1},u_n)\le t, t_{n-1}\ge t_n\}
\mu_{n+1}(d(u_0,\ldots,u_n))\notag\\
\le& 2(t^d-s^d)\iint f_n V_d(K_{n})\mu_{n}(d(u_0,\ldots,u_{n-1}))
\BQ(d(t_{n},K_{n}))\notag\\
\label{sum2}
&+2\iint f_n h(s,t)\sum_{i=0}^{d-1}\kappa_d\binom{d}{i}(t_{n-1}-t_n)^i \rho(K_{n})^i
(\rho(K_{n-1})+\rho(K_{n-1}))^{d-i}\notag \\
&\qquad\qquad \mu_{n}(d(u_0,\ldots,u_{n-1}))\BQ(d(t_{n},K_{n})).
\end{align}
Combining \eqref{sum11}-\eqref{sum2} with \eqref{sum0} yields
\begin{align*}
&\BP(\Psi\in F_n(s,t,B))\le c_1^{n+1}h(s,t)c_2\iint f_n((t_{n-1}\vee 1)+(t_n\vee 1))^d \\
&\qquad\qquad
((\rho(K_{n-1})\vee 1)+(\rho(K_n)\vee 1))^{d}\mu_{n}(d(u_0,\ldots,u_{n-1}))\BQ(d(t_{n},K_{n})).
\end{align*}
Since the integrand in \eqref{sum0} equals $f_{n+1}$, we obtain recursively
\begin{align}\label{aaa3}
&\BP(\Psi\in F_n(s,t,B)) \le V_d(B) c_3^nc_4h(s,t)^n\\
&\int \prod_{i=1}^{n}((t_{i-1}\vee 1)+(t_i\vee 1))^d 
((\rho(K_{i-1})\vee 1)+(\rho(K_i)\vee1))^{d}
\BQ^{n+1}(d((t_0,K_0),\ldots,(t_n,K_n))).\notag
\end{align}

It is convenient to introduce a random vector $((\tau_0',Y_0),\ldots,(\tau_n',Y_n))$
with distribution $\BQ^{n+1}$ and to define $\tau_i:=\tau_i'\vee 1$ and
$\rho_i:=\rho(Y_i)\vee 1$ for $i=0,\ldots,n$. Then inequality \eqref{aaa3} can be 
written as
\begin{align*}
\BP(\Psi\in F_n(s,t,B))
\le  V_d(B)c_3^nc_4h(s,t)^n \BE \Big[\Big(\prod_{i=1}^{n}(\tau_{i-1}+\tau_i)\Big)^d
\Big(\prod_{i=1}^{n}(\rho_{i-1}+\rho_i)\Big)^{d}\Big].
\end{align*}
Now we can use the elementary inequality
$$
\prod^n_{i=1} (a_{i-1}+a_i)\le 2^{n}\Big(\prod^{n}_{i=1} a_{i-1}a_{i}\Big)
\le 2^{n}\prod^{n}_{i=0}a^2_i
$$
that holds whenever $a_0,\ldots,a_n\ge 1$. This gives
\begin{align*}
\BP(\Psi\in F_n(s,t,B))
&\le  V_d(B)c_3^nc_4h(s,t)^n 2^{2nd}\BE \Big[\prod_{i=0}^{n}\tau_{i}^{2d} \rho_i^{2d}\Big]\\
&=V_d(B)c_3^nc_4h(s,t)^n (4^d)^n(\BE[\tau_{0}^{2d} \rho_0^{2d}])^{n+1}.
\end{align*}
Assumption \eqref{assmoment} implies that
\begin{align*}
\BP(\Psi\in F_n(s,t,B))\le V_d(B) c_5 c_{6}^nh(s,t)^n.
\end{align*}
Now we define $b_i:=(i/(c_6+1))^{1/d}$, $i\in\N_0$. By definition \eqref{defh} of $h$
\begin{align*}
h(b_{i-1},b_i)\le \frac{1}{c_{6}+1}.
\end{align*}
Hence $\lim_{n\to\infty}\BP(\Psi\in F_n(b_{i-1},b_i,W_j))=0$, $i,j\in\N$,
so that \eqref{defF1} implies $\BP(\Psi\in F)=0$.

We now show that \eqref{c2}, \eqref{c1} and \eqref{c3} have full measure.
The probability that $\Psi$ does not satisfy \eqref{c2} can be bounded by
\begin{align}\label{aaa6}
\BE\sum_{(u_0,u_1,u_2)\in\Psi^{(3)}}\I\{d(u_0,u_1)=d(u_0,u_2)\}
+\BE\sum_{(u_0,u_1,u_2,u_2)\in\Psi^{(4)}}\I\{d(u_0,u_1)=d(u_2,u_3)\},
\end{align}
where the sum is taken over all triples resp.\ quadruples 
of points of $\Psi$ with pairwise different entries. We recall the definition 
\eqref{deffacmommeas} of the factorial moment measure $\nu^{(k)}$ for $k\in\N$  
and that \eqref{assfac} holds for $\Psi$.
Together with the definition \eqref{d}  of  $d$
it follows that the first summand of \eqref{aaa6} is bounded by
\begin{align*}
a \iiint \I\{x_2\in \partial A(u_0,u_1,(s_2,K_2))\} dx_2\BQ(d(s_2,K_2))\mu_2(d(u_0,u_1)),
\end{align*}
where $A(u_0,u_1,(s_2,K_2))$ is a convex body, depending on $u_0,u_1$
and $(s_2,K_2)$ but not on $x_2$. Hence this integral is zero.
The second summand in \eqref{aaa6} can be treated analogously.
Thus we obtain that $\Psi$ satisfies \eqref{c2} almost surely.
That $\Psi$ does also satisfy \eqref{c3} almost surely
can be shown with the same type of argument.

Finally we prove that
$\BP(\card\{v\in\Psi\setminus\{u\}:d(u,v)\le r\}<\infty \mbox{ for all } u\in\Psi)=1$
for all $r>0$. It is sufficient to show, that
\begin{align}\label{aaa5}
\BE\sum_{(u,v)\in\Psi^{(2)}}\I\{\card\{v\in\Psi\setminus\{u\}:d(u,v)\le r\}<\infty\}
\I\{u\in W_j\times\R_+\times\mathcal{K}^d\}
\end{align}
is finite for all $j\in\N$ and $r>0$. For fixed $j$ and $r$ this expectation
is bounded by
\begin{align*}
a^2\iiint \I\{x\in W_j,d((x,s,K),(y,t,L)\le r\} dxdy\BQ^2(d((s,K),(t,L))),
\end{align*}
which in turn is bounded by
\begin{align}\label{aaa4}
a^2&\iiint \I\{x\in W_j,s\le t, s+a(x,K,y,0)\le r\} dxdy\BQ^2(d((s,K),(t,L)))\notag\\
+&a^2\iiint \I\{x\in W_j,s\le t, t+a(x+(t-s)K,K,y,L)\le r\} dxdy\BQ^2(d((s,K),(t,L)))\notag\\
+&a^2\iiint \I\{x\in W_j,t\le s, t+a(y,L,x,0)\le r\} dxdy\BQ^2(d((s,K),(t,L)))\notag\\
+&a^2\iiint \I\{x\in W_j,t\le s, s+a(y+(s-t)L,L,x,K)\le r\} dxdy\BQ^2(d((s,K),(t,L))).
\end{align}
The first two summands equal
\begin{align*}
a^2&\iiint \I\{x\in W_j,s\le t, s\le r\} \I\{y\in x+(r-s)K\} dxdy\BQ^2(d((s,K),(t,L)))\\
+&a^2\iiint \I\{x\in W_j,s\le t, t\le r\} \I\{y\in x+(t-s)K+(r-t)(K-L)\}\\
&\qquad\qquad\qquad\qquad\qquad\qquad\qquad\qquad\qquad\qquad\qquad\qquad\qquad 
dxdy\BQ^2(d((s,K),(t,L))).
\end{align*}
This sum is bounded by
\begin{align*}
a^2&\iiint \I\{x\in W_j,s\le t, s\le r\} \kappa_d (r-s)^d\rho(K)^d dxdy\BQ^2(d((s,K),(t,L)))\\
+&a^2\iiint \I\{x\in W_j,s\le t, t\le r\}\kappa_d [(t-s)\rho(K)+(r-t)(\rho(K)+\rho(L))]^d\\
&\qquad\qquad\qquad\qquad\qquad\qquad\qquad\qquad\qquad\qquad\qquad\qquad\qquad 
dxdy\BQ^2(d((s,K),(t,L))).
\end{align*}
By \eqref{assmoment} both summands are finite. By symmetry the third and 
fourth summand in \eqref{aaa4} can be treated in the same way. 
Therefore \eqref{aaa5} is finite for all $j\in\N$ and $r>0$.\qed

\section{Absence of percolation}\label{secperc}

In this section we consider a point process $\Psi=\{(X_n,T_n,Z_n):n\ge 1\}$
satisfying the moment assumptions \eqref{assfac}, where the
probability measure $\BQ$ is assumed to satisfy \eqref{assmoment} and,
moreover, to be concentrated on the set 
of all strictly convex bodies.
In addition we assume $\Psi$ is non-empty and {\em stationary}, that is 
$\BP(x+\Psi\in\cdot)$ does not depend on $x\in\R^d$.
The {\em intensity} of $\Psi$ is defined by
$\gamma_\Psi:=\BE \Psi([0,1]^d\times\R_+\times \cK)$. By \eqref{assfac}
this is a finite number, while $\Psi\ne\emptyset$ implies $\gamma_\Psi>0$.

Due to Theorem \ref{texist} there exists a $\BP$-a.s.\ unique
growth-maximal hard-core model $\Psi^*=\{(X_n,T_n,Z_n,R_n):n\ge 1\}$
based on $\Psi$. We define
\begin{align}\label{Z}
Z:=Z(\Psi):=\bigcup_{n=1}^\infty X_n+R_nZ_n=\bigcup_{n: R_n>0} X_n+R_nZ_n,
\end{align}
that is the union of all grains (which started growing). Note that $X_n\in Z$ if $R_n=0$. 
We say that
$Z$ (or the growth-maximal hard-core model) {\em percolates},
if $Z$ contains an unbounded connected component.
To be more exact we introduce a graph with vertex set
$\Psi_+:=\{(X_n,T_n,Z_n): R_n>0\}$. Two different points 
$(X_m,T_m,Z_m), (X_n,T_n,Z_n)\in\Psi_+$ share an edge if
they are neighbours, that
is if  $(X_m+R_mZ_m)\cap(X_n+R_nZ_n)\ne \emptyset$.
A {\em cluster} is a connected component of this graph.

As $\BP(\Psi\neq\emptyset)=1$, the second defining property (ii)
of a growth-maximal hard-core model given in the introduction together
with Theorem \ref{texist} implies that almost surely $\Psi_+$ is non-empty too.
Since $\Psi_+$ is stationary (by the translation invariance of $R$)
we have in fact that $\BP(\card\Psi_+=\infty)=1$.
Our aim in this section is to verify the following theorem.

\begin{theorem}\label{tperc} Almost surely there are no infinite clusters.
\end{theorem}

This theorem implies in particular that the random set $Z$
does not percolate. Our proof is based on some ideas in \cite{HM96}
and \cite{DL05}. We begin with the following lemma.

\begin{lemma}\label{learlier} Almost surely any point
$(X_n,T_n,Z_n)\in\Psi_+$ has exactly one
earlier neighbour.
\end{lemma}
{\sc Proof:} Assume without loss that $\Psi\in H$ with $H$ being the set
from Proposition \ref{p2.2}. In addition suppose $u=(x,s,K)\in\Psi$ with
$R(u):=R(\Psi,u)>0$ where $R(\cdot,\cdot)$ is the function introduced in
Proposition \ref{p2.2}. By growth-maximality $u$ has at least
one earlier neighbour.
If $u$ has more than one earlier neighbour, then for two different neighbours
$v=(y,t,L)$ and $w=(z,r,M)$ of $u$ one of the following cases must occur:
\begin{align}
R(u)&=R(v)+t-s=R(w)+r-s>0;\label{ca1}\\
R(u)&=R(v)+t-s>0,\ R(u)+s>R(w)+r>0;\label{ca2}\\
R(u)&+s>R(v)+t>0,\ R(u)+t>R(w)+r>0\label{ca3}.
\end{align}
We need to show that all three cases \eqref{ca1}-\eqref{ca3} have probability zero. 
In order to understand the argument we first illustrate each of the three cases 
under the additional assumption $t\le s\le r$. Because of the construction 
algorithm from the proof of Proposition \ref{p2.2} it follows from \eqref{ca1} that
$R(u)=a(x+(t-s)K,K,y,L)=a(x,K,z+(r-s)M,M)+r-s$.
In the case \eqref{ca2} we have $R(u)=a(x+(t-s)K,K,y,L)=a(x,K,z+R(w)M,0)$. 
Moreover there exists an $n\in\N$ such that $R(w)$ can be replaced
by an expression depending on $v$ and a neighbour 
$w_1=(z_1,r_1,M_1)\in\Psi_n\setminus\{u,v\}$, 
where $\Psi_n$ is defined in the construction algorithm in the 
proof of Proposition \ref{p2.2}. For $s\ge r_1$ we have either
$R(w)=a(x+(s-r_1)M,M,z_1,M_1)$ or $R(w)=a(z,M,z_1+R(w_1)M_1,0)$. If $s<r_1$ 
we have either that $R(w)=a(x,M,z_1+(r_1-s)M_1,M_1)+(r_1-s)$ or 
that $R(w)=a(z,M,z_1+R(w_1)M_1,0)$. 
Now $R(w_1)$ can be replaced in the same way using $w_1$ and a neighbour 
$w_2\in\Psi_m\setminus\{u,v,w_1\}$, $m\le n$. 
After a finite number of steps this procedure ends 
and we have $a(x,K,z+R(w)M,0)=a(x,K,z+f(w,w_1,\ldots,w_k)M,0)$ for some $k\ge 0$.
The third case leads to $R(u)=a(x,K,y+R(v)L,0)=a(x,K,z+R(w)M,0)$. 
We can apply the same replacing routine as in the second case for $R(v)$ and $R(w)$.

In the case \eqref{ca1} we have that
\begin{align}\label{gf3}
a(x+(t-s)_+K,&K,y+(t-s)_+L,L)\\ 
\notag
&=a(x+(r-s)_+K,K,z+(s-r)_+M,M)+(r-s)_+-(t-s)_+>0,
\end{align}
where $b_+:=\max\{b,0\}$ for $b\in\R$.
In the other two cases the preceding replacement process yields countable 
families $F_n, n\ge 1,$ of functions $f:(\R^d\times\R_+\times\cK)^n\rightarrow\R_+$ 
such that the following is true. There exist $k,l,m\ge 1$, 
$f_1\in F_k,f_2\in F_l,f_3\in F_{m}$ and 
$(v_1,\ldots,v_k)\in(\Psi\setminus\{u\})^{(k)}$, 
$(w_1\ldots,w_l)\in(\Psi\setminus\{u\})^{(l)}$ as well as 
$(w'_1,\ldots,w'_{m})\in(\Psi\setminus\{u\})^{(m)}$ such that
\begin{align}\label{gf2}
  a(x+(t-s)_+K,K,y+(t-s)_+L,L)+(t-s)_+
=a(x,K,z+f_1(v_1,\ldots,v_{k})M,0)>0.
\end{align}
or
\begin{align}\label{gf}
a(x,K,y+f_2(w_1,\ldots,w_{l})L,0)
=a(x,K,z+f_3(w'_1,\ldots,w'_{m})M,0)>0,
\end{align}

Using that the set $F:=\cup_{n\in\N}F_n$ is countable and assumption \eqref{assfac} 
on the factorial moment measures it is not difficult to see that the 
probabilities of all three cases \eqref{gf3}-\eqref{gf2} are zero, provided 
the set of $x\in\R^d$ satisfying \eqref{gf3}-\eqref{gf} is of Lebesgue measure zero. 
In the case \eqref{gf} this means that
\begin{align}\label{gf4}
\int\I\{a(x,K,K',0)=a(x,K,K'',0)>0\} dx =0,
\end{align}
whenever $K$ is strictly convex and the interiors
of $K'$ and $K''$ do not intersect. In fact, if
the previous properties are satisfied and
$r:=a(x,K,K',0)=a(x,K,K'',0)>0$ for some $x\in\R^d$,
then the infimum $r=\inf\{t>0:(x+tK)\cap (K'\cup K'')\ne\emptyset\}$
is attained in two different points in the boundary
of $K'\cup K''$. It follows from the Lipschitz property
of the $K$-distance and the strict convexity of $K$
that the set of points $x$ with this property
has Lebesgue measure $0$. 
The detailed argument (even for sets more general than $K'\cup K''$)
can be found in \cite{Hug99}.
The cases \eqref{gf3} and \eqref{gf2} can be treated similarly.
\qed

\bigskip

Two different points $(X_m,T_m,Z_m),(X_n,T_n,Z_n)\in\Psi_+$ form a {\em doublet}
if 
$T_m+R_m=T_n+R_n$. This means that the points stop each other mutually.

\begin{lemma}\label{ldoublet} Almost surely any cluster
contains at most one doublet while any finite cluster
contains exactly one doublet.
\end{lemma}
{\sc Proof:} For $u=(x,s,K)\in\Psi$ we
write $R(u):=R(\Psi,u)$ and $S(u):=s+R(u)$.
Consider a cluster $\Xi\subset\Psi_+$ and assume that
$\{u,v\},\{u',v'\}$ are two different doublets in $\Xi$.
Lemma \ref{learlier} implies that $\{u,v\}\cap\{u',v'\}=\emptyset$.
(In the following we ignore $\BP$-null sets.)
Assume (without loss) that $S(u)\ge S(u')$ and
let $u_0,u_1,\dots,u_m$ be a path in $\Psi_+$ such that
$u_0=u$, $u_1\notin\{u,v\}$, $u_m=u'$, and $v'\notin\{u_1,\ldots,u_m\}$.
Lemma \ref{learlier} implies that $S(u)<S(u_1)$ (otherwise
$u$ would have two earlier neighbours $v$ and $u_1$)
and then, recursively,
$S(u_i) < S(u_{i+1})$ for all $i\in\{0,\ldots,m-1\}$.
Therefore $S(u_{m-1})<S(u_m)=S(u')=S(v')$, contradicting the fact
that $u'$ has only one earlier neighbour. Therefore
the two doublets cannot be connected by a path,
so that $\Xi$ can have at most one doublet. 
If $\Xi$ is finite, then
any $u\in\Xi$ that minimizes the function $S$ on $\Xi$ must
be a member of a doublet. \qed

\bigskip
{\sc Proof of Theorem \ref{tperc}:}
Assume that $\Xi\subset\Psi_+$ is an infinite cluster without a
doublet. Take $u_1\in \Xi$ and let $u_2$ be the (unique) earlier
grain neighbour of $u_1$. By assumption we have $S(u_1)>S(u_2)$.
Continuing this way we obtain an infinite sequence $u_1,u_2,\dots$
of different points in $\Psi$ such that $R(u_n)>0$ and
$S(u_n)>S(u_{n+1})$ for all $n\ge 1$. Now assume that
there are $s<t$ such that $s< S(u_n)\le t$ for all $n\ge 1$.
Then Lemma \ref{l2.1} implies that $s< d(u_n,u_{n+1})\le t$
for all $n\ge 1$. We have shown in the proof of Theorem \ref{texist}
that his event has probability zero.

In view of Lemma \ref{ldoublet} it remains to prove that
there are no infinite clusters with exactly one doublet.
It is here, where stationarity plays a crucial role.
In contrast to \cite{HM96,DL05} we use the
{\em mass-transport principle} (see e.g.\cite{Last})
\begin{align}\label{mtp}
\BE \iint \I\{u\in B\}g(\Psi_+,u,v)\Psi_+(dv)\Psi_+(du)
=\BE \iint \I\{v\in B\}g(\Psi_+,u,v)\Psi_+(du)\Psi_+(dv),
\end{align}
that holds for all measurable
$B\subset \R^d\times\R_+\times\cK$ and all
non-negative measurable functions $g$,
provided that $g$ is translation-invariant.
For $v\in \Psi_+$ define $C(v)\equiv C(\Psi_+,v)\subset\Psi_+$
as the cluster containing $v$.
Let $\Psi_\infty$ be the set of all $v\in\Psi_+$ such that
$C(v)$ is infinite and has a doublet.
Define $g(\Psi_+,u,v):=1$ if $u\in\Psi_\infty$, $v\in C(u)$,
and if $u$ is the lexicographically smallest point
of the doublet of $C(v)$. Otherwise, set $g(\Psi_+,u,v):=0$.
With this choice of $g$ and $B:=[0,1]^d\times\R_+\times\cK$
the right-hand side of
\eqref{mtp} is at most $\BE\Psi_+[0,1]^d\le\gamma_\Psi$.
The left-hand side vanishes if the intensity $\gamma_\infty$ of
the stationary point process $\Psi_\infty$ is zero.
Otherwise $\Psi_\infty$ has infinitely many points and
the left-hand side is infinite. This shows that
$\gamma_\infty=0$, as asserted.\qed

\section{A central limit theorem}\label{secstabclt}

In this section we assume that $\Psi$ is an independently
marked Poisson process on $\R^d$ with intensity $1$ and
mark space $[0,\infty)\times\cK$
satisfying  $\BP$-a.s.\ the restrictions $T_0=0$ and
\begin{align}\label{shapecond}
B^d\subset Z_0\subset cB^d
\end{align}
for some fixed $c\ge 1$. The assumption $T_0=0$ means that
all grains are born at the same time (taken as $0$ without
further restriction of generaliy), while
\eqref{shapecond} implies that all growth times are strictly positive and
bounded from below and above independently of the grain shape.
Using stabilization arguments (\cite{LaPe12,PeEJP})
we shall prove a central limit theorem for the 
growth times $R_n$, $n\in\N$.

Let $g$ be a finite kernel from $\R_+\times\R^d\times\cK_0$ to $\R^d$,
where $\cK_0$ is the set of all $K\in\cK$ satisfying $B^d\subset K\subset c B^d$.
We assume $g$ to be translation invariant, that is
$g(t,x,K,A)=g(t,x+y,K,A+y)$ for all $y\in\R^d$. In addition we let $\bN_0$
denote the set of all $\psi\in\bN$, such that for all $(x,s,K)\in\psi$ we have $s=0$
and $K\in\cK_0$. We use the abbreviation $(x,K):=(x,0,K)$ for all $(x,0,K)\in\psi$
and $\psi\in\bN_0$ and suppress the birth time in all other expressions as well.
For all $x,y\in\R^d$ and $K,L\in\cK_0$ define $\psi^{(x,K)}:=\psi\cup\{(x,K)\}$
and $\psi^{(x,K),(y,L)}:=\psi\cup\{(x,K),(y,L)\}$ if such a union lies in $\bN_0$
and $\psi^{(x,K)}:=\psi$ resp.\ $\psi^{(x,K),(y,L)}:=\psi$ otherwise. Set
\begin{align}\label{defrho}
\rho(\psi,x,K,A):=g(R(\psi^{(x,K)},x,K),x,K,A)\I_{\{\card(\psi^{(x,K)})> 1\}}
\end{align}
for all $\psi\in\bN_0$, $x\in\R^d$, $K\in\cK_0$ and measurable $A\subset\R^d$,
where $R$ is the function defined in Proposition \ref{p2.2}.
With
\begin{align}
\mu(\psi,A):=\sum_{(x,K)\in\psi} \rho(\psi,x,K,A),
\end{align}
for $\psi\in\bN_0$ and measurable $A\subset\R^d$ we define a measure 
$\mu(\psi)(\cdot):=\mu(\psi,\cdot)$. Define the {\em observation windows} 
$W_1:=[-1/2,1/2]^d$ and 
$W_n:=n^{1/d}W_1$ for $n\ge 2$. Furthermore set $\psi_A:=\psi\cap(A\times\cK_0)$.

\begin{theorem}\label{clt}
Assume there exist $\alpha,\beta>0$ such that the translation invariant kernel 
$g$ in \eqref{defrho} satisfies the growth bound
\begin{align}\label{gb}
g(t,x,K,\R^d)\le\alpha t^\beta, \quad t\ge 0,\ x\in\R^d,\ K\in\cK_0.
\end{align}
Moreover suppose that $f:W_1\to \R$ is a bounded, almost everywhere continuous function. 
Then the limit
$ \sigma_{\mu,f}:=\lim_{n\to\infty}\BV(\int_{W_1}f d\mu_n^\rho)/n$
exists and as $n\to\infty$
\begin{align}\label{cltstatement}
\frac{1}{\sqrt{n}}\Big(\int_{W_n} f(n^{-1/d}x)\mu(\Psi_{W_n},dx)
-\BE\Big[\int_{W_n} f(n^{-1/d}x)\mu(\Psi_{W_n},dx)\Big]\Big)\stackrel{d}{\to}
\mathcal N(0,\sigma_{\mu,f}).
\end{align}
\end{theorem}
We will derive this theorem from Theorem 2.1 and 2.2 from \cite{PeEJP}
starting with showing that the growth times 
$R(\psi^{(y,L)},y,L)$ stabilize. Here we borrow heavily from \cite{LaPe12} several times.

For $\psi\in\bN_0$, $y\in\R^d$ and $K\in\cK_0$ define
\begin{align}\label{defD}
D(\psi,y,L):=\inf\{a(y,L,x,0): (x,K)\in\psi\setminus\{(y,L)\}\},
\end{align}
with $a$ as in \eqref{a}. If $(y,L)\in\psi$ then $D(\psi,y,L)$ 
is an upper bound for the growth time $R(\psi,y,L)$ of the grain $(y,L)$.

We recall the function $d$ defined in \eqref{d} and the constant $c\ge 1$ 
from assumption \eqref{shapecond}. Below we will use that
\begin{align}\label{propd}
\frac{1}{2c}\|x-y\|\le d(x,K,y,L)\le \frac{1}{2}\|x-y\|,\quad x,y\in\R^d,\ K,L\in\cK_0,
\end{align}
provided that $B^d\subset K\subset cB^d$ and $B^d\subset L\subset cB^d$.

By a \textit{(finite) descending chain} in $\psi\in\bN_0$ we mean a 
finite sequence $u_0,\ldots,u_n$ $(n\ge 1)$ of distinct points of $\psi$ for 
which $d(u_{i-1},u_i)\ge d(u_i,u_{i+1})$ for all $i\in\{1,\ldots,n-1\}$. 
For $\psi\in\bN_0$, $y\in\R^d$ and $L\in\cK_0$ a triple 
$(x,K,r)\in\R^d\times\cK_0\times(0,\infty)$ belongs to the set $A(\psi,y,L)$ 
if there exists a descending chain $u_0,\ldots,u_n$ in $\psi^{(y,L)}$, 
such that $u_0=(y,L)$, $u_n=(x,K)$, $d(u_0,u_1)\le D(\psi^{(y,L)},y,L)$ and 
$r=d(u_{n-1},u_n)$. We define
\begin{align*}
S(\psi,y,L):=B(y,2cD(\psi^{(y,L)},y,L))\cup\bigcup_{(x,K,r)\in A(\psi,y,L)} B(x,2cr),
\end{align*}
if $\psi\setminus\{(y,L)\}\neq\emptyset$ and $S(\psi,y,L)=\R^d$ otherwise and 
$S^*(\psi,y,L):=S(\psi,y,L)\times\cK_0$.

\begin{lemma}\label{l4.3}
Let $\psi,\varphi\in\bN_0$, $y\in\R^d$ and $L\in\cK_0$. If 
$\psi\cap S^*(\psi,y,L)=\varphi\cap S^*(\psi,y,L)$ then $S(\psi,y,L)=S(\varphi,y,L)$.
\end{lemma}
{\sc Proof:} We assume $\psi\setminus\{(y,L)\}\neq\emptyset$ since otherwise 
the result is trivial. Since $B((y,2cD(\psi^{(y,L)},y,L))\subset S(\psi,y,L)$ 
we have $D(\psi^{(y,L)},y,L)=D(\varphi^{(y,L)},y,L)$.

Take $(x,K,r)\in A(\psi,y,L)$. Thus there exists a descending chain 
$u_0,\ldots,u_n\in\psi^{(y,L)}$ such that $u_0=(y,L)$, $u_n=(x,K)$, 
$d(u_0,u_1)\le D(\psi^{(y,L)},y,L)$ and $r=d(u_{n-1},u_n)$. We have 
$u_m\in S^*(\psi^{(y,L)},y,L)$ for $1\le m\le n$. Therefore $u_m\in\varphi$. 
So $(x,K,r)\in A(\varphi,y,L)$ and hence $A(\psi,y,L)\subset A(\varphi,y,L)$.

Consider $(x,K,r)\in A(\varphi,y,L)$. So there exists a descending chain 
$u_0,\ldots,u_n\in\varphi^{(y,L)}$ with
$u_0=(y,L)$, $u_n=(x,K)$, $d(u_0,u_1)\le D(\varphi^{(y,L)},y,L)$ and $r=d(u_{n-1},u_n)$. 
We show $u_m\in\psi$ for $1\le m\le n$ by induction on $m$. For $m=1$ we have 
$d(u_0,u_1)\le  D(\varphi^{(y,L)},y,L)$. Using \eqref{propd} with $(y,L)$ and 
$u_1=(x_1,K_1)$ we obtain $\|y-x_1\|\le 2cD(\varphi^{(y,L)},y,L)$. Therefore 
$u_1\in \varphi\cap S^*(\psi,y,L)$ and hence $u_1\in\psi$. Now assume our 
assertion holds for $u_m$, $1\le m\le k-1$. Applying \eqref{propd} with 
$u_m=(x_m,K_m)$ and $u_{m+1}=(x_{m+1},K_{m+1})$ yields $\|x_m-x_{m+1}\|\le 2cd(u_m,u_{m+1})$. 
So $u_k\in \varphi\cap S^*(\psi,y,L)$ and thus $u_k\in\psi$ 
so that $A(\varphi,y,L)\subset A(\psi,y,L)$. Together with the inclusion 
proved above we deduce $A(\psi,y,L)= A(\varphi,y,L)$ and the assertion follows.\qed

\bigskip
We recall from the proof of Lemma \ref{l2.3} the definition of a measurable set 
$G\subset H$ such that $\BP(\Psi\in G)=1$.
\begin{lemma}\label{l4.1}
Suppose $\psi\in\bN_0\cap G$ and $(y,L)\in\R^d\times\cK_0$. If $S(\psi,y,L)$ 
is bounded, then it satisfies
\begin{align*}
R(\psi^{(y,L)},y,L)=R(\psi^{(y,L)}\cap S^*(\psi,y,L)\cup\varphi,y,L),
\end{align*}
for all finite $\varphi\in\bN_0$ with $\varphi\subset\R^d\setminus S(\psi,y,L)\times\cK_0$.
\end{lemma}
{\sc Proof:} Suppose $S(\psi,y,L)$ is bounded and $\psi\in\bN_0$. 
Since $S(\psi^{(y,L)},y,L)=S(\psi,y,L)$ we assume w.l.o.g. $(y,L)\in\psi$. 
Moreover because of Lemma \ref{l4.3} it suffices to prove the equality in the case 
$\psi\subset S^*(\psi,y,L)$. Hence we assume this too. We note that since 
$\varphi$ is finite the set $\psi':=\psi\cup\varphi$ is an element of $\bN_0\cap G$ 
and define $R(x,K):=R(\psi,x,K)$ and $R'(x,K):=R(\psi',x,K)$.

Suppose $R(y,L)>R'(y,L)$ and that $(x_1,K_1)$ is an earlier neighbour of $(y,L)$ 
in $\psi'$. We assume for now $(x_1,K_1)\in\psi$. Since $(x_1,K_1)$ and 
$(y,L)$ are neighbours we obtain $(y+R'(y,L)L)\cap(x_1+R'(x_1,K_1)K_1)\neq\emptyset$. 
Because of the hard-core property on $\psi$ and $R(y,L)>R'(y,L)$ it follows 
$R(x_1,K_1)<R'(x_1,K_1)$. Suppose $(x_2,K_2)$ is an earlier neighbour of $(x_1,K_1)$ 
in $\psi$. The hard-core property on $\psi'$ yields $R(x_2,K_2)>R'(x_2,K_2)$. 
Let $(x_3,K_3)$ be an earlier neighbour of $(x_2,K_2)$ in $\psi'$ and assume 
$(x_3,K_3)\in\psi$. Continuing this procedure leads to a sequence of points satisfying
\begin{align}\label{RandR'}
R(y,L)>R'(y,L)\ge R'(x_1,K_1)>R(x_1,K_1) \ge R(x_2,K_2)>R'(x_2,K_2)>R'(x_3,K_3)\ldots,
\end{align}
terminating at $(x_n,K_n)$ if $(x_n,K_n)\in\varphi$ (so $n$ must be odd). 
The strict inequalities ensure that the points of this (possibly) terminating 
sequence are all different. By Lemma \ref{l2.1} we get
\begin{align}
R(x_i,K_i)\ge d((x_i,K_i),(x_{i+1},K_{i+1}))\ge R(x_{i+1},K_{i+1})
\end{align}
for odd $i\in\N$ and the same inequalities if we replace $R$ by $R'$ and $i\in\N$ is even. 
Furthermore $R'(y,L)\ge d((y,L),(x_1,K_1))\ge R'(x_1,K_1)$ holds. 
Therefore \eqref{RandR'} implies
\begin{align*}
d((y,L),(x_1,K_1))\ge d((x_1,K_1),(x_2,K_2))\ge d((x_2,K_2),(x_3,K_3))\ge\ldots.
\end{align*}
Thus the sequence forms a descending chain with $d((y,L),(x_1,K_1))\le D(\psi^{(y,L)},y,L)$. 
If the sequence terminates at some point $(x_n,K_n)\in\psi'$ 
then $(x_n,K_n)\in S^*(\psi,y,L)$ contradicting $\varphi\cap S^*(\psi,y,L)=\emptyset$. 
If the sequence does not terminate then $S(\psi,y,L)$ cannot be bounded which is 
also a contradiction. Hence if $R(y,L)>R'(y,L)$ the assertion is shown. The case 
$R'(y,L)>R(y,L)$ is proven similarly, this time starting with an earlier 
neighbour $(x_1,K_1)$ of $(y,L)$ in $\psi$.\qed

\bigskip

Let $B(y,r):=\{x\in\R^d:\|x-y\|\le r\}$ be the ball around $y\in\R^d$ of radius $r> 0$ 
and $B^*(y,r):=B(y,r)\times\cK_0$.
Furthermore define for $\psi\in\bN_0\cap G$, $y\in\R^d$ and $K\in\cK_0$
\begin{align}\label{defU}
U(\psi,y,L):=\inf\{r> 0: S(\psi,y,L)\subset B(y,r)\}
\end{align}
and $U(\psi,y,L):=\infty$ otherwise. Because of the Lemma \ref{l4.1} we have
\begin{align}\label{stabradforR}
R(\psi^{(y,L)},y,L)=R(\psi^{(y,L)}\cap B^*(y,U(\psi,y,L)\cup\varphi,y,L),
\end{align}
for all finite $\varphi\in\bN_0$ with 
$\varphi\subset\R^d\setminus B(y,U(\psi,y,L)\times\cK_0$. Hence we call $U$ a 
{\em radius of stabilization of $R$ at $(y,L)$ with respect to $\psi$}. 
In \cite{PeEJP} the radius of stabilization $T$ of $\rho$ with respect to 
$\psi$ and $K$ is defined as
\begin{align}\notag
T(\psi,K):=\inf\{r\ge 0:&\rho([\psi\cap B^*(0,t)\cup\varphi],0,K,\R^d)
=\rho([\psi\cap B^*(0,t)],0,K,\R^d)\\\label{stabrad}
&\ \ \qquad\mbox{ for all finite } \varphi\in\bN_0 \mbox{ with } 
\varphi\subset(\R^d\setminus B(0,t))\times\cK_0\}.
\end{align}
By definition $\rho(\psi,0,K,\R^d)$ depends on $\psi$ only via the radius 
$R(\psi^{(0,K)},0,K)$. Because of \eqref{stabradforR} it follows that $T(\psi,K)$ 
is bounded by $U(\psi,0,K)$. Now we deal with the tail behaviour of $U(\Psi,0,Z_0)$.
\begin{lemma}\label{l4.2}
Let $0<\gamma<1$. Then there exist constants $\alpha,\beta>0$ such that
\begin{align*}
\BP(U(\Psi,0,Z_0)>t)\le \alpha\exp(-\beta t^\gamma),\quad t\ge 0.
\end{align*}
The same holds for $U(\Psi^{(x,Z_1)},0,Z_0)$ for all $x\in\R^d$, where $Z_1$ 
is an independent copy of $Z_0$ and also independent of $\Psi$.
\end{lemma}
{\sc Proof:} We start with the first assertion.  W.l.o.g we assume $t\ge 1$. 
Define $U:=U(\Psi,0,Z_0)$ and $D:=D(\Psi,0,Z_0)$. Let $0<\varepsilon<1/(d+1)$. It holds
\begin{align*}
\{U>t\}\subset 
\Big\{D >\frac{t^\varepsilon}{4c}\Big\}\cup\Big\{D\le \frac{t^\varepsilon}{4c},U>t\Big\},
\end{align*}
where $c>1$ is the constant from the shape condition \eqref{shapecond}. Define $E$ 
as the set of all $\psi\in\bN_0$, which contain a descending chain 
$(0,K_0),u_1,\ldots,u_n$ such that
\begin{align}\label{dcond}
d((0,K_0),u_1)\le D(\psi,0,K_0)\le t^\varepsilon/(4c)
\end{align}
and $n\ge t^{1-\varepsilon}/2$. For $\psi\in\bN_0$ and $K_0\in\cK_0$ such that 
$\psi^{(0,K_0)}\in\bN_0$, $U(\psi,0,K_0)>t$ and \eqref{dcond} holds
we have $\psi^{(0,K_0)}\in E$. To see this, assume on the contrary that every 
descending chain of $\psi^{(0,K_0)}$ which satisfies \eqref{dcond} would 
consist of $n$ points with $n<t^{1-\varepsilon}/2$. Since
\begin{align*}
S(\psi,0,K_0)\subset B(0,2c(n+1)D(\psi,0,K_0))
\end{align*}
this would imply
\begin{align*}
U(\psi,0,K_0)\le 2c\Big(\frac{t^{1-\varepsilon}}{2}+1\Big)D(\psi,0,K_0)
\end{align*}
by the definition \eqref{defU} of $U$. Since $D(\psi,0,K_0)\le t^\varepsilon/(4c)$, 
$t\ge 1$, and $0<\varepsilon<1/(d+1)$ this would lead to the contradiction 
$U(\psi,0,K_0)\le (3/4)t$. So we deduce
\begin{align}\label{ineq4.3}
\BP(U>t)\le \BP\Big(D >\frac{t^\varepsilon}{4c}\Big)+\BP\Big(\Psi^{(0,Z_0)}\in E\Big).
\end{align}
Define $k:=\lceil(c_1+1)t^{\varepsilon d}/(4c)^d\rceil$ and $s_i:=(i/(c_1+1))^{1/d}$, 
$i\ge 0$, where $c_1:=V_d(B(0,2c))$ and $\lceil\cdot\rceil$ is the ceiling function. 
Hence $s_k\ge t^\varepsilon/(4c)$. Set $m:=\lceil t^{1-\varepsilon}/(2k)\rceil$ and 
$n_0:=\lceil t^{1-\varepsilon}/2\rceil$. Note that if $t\to\infty$ then $m\to\infty$ 
as well. W.l.o.g assume $t\ge 1$ so large that $m\ge 3$. For $\psi\in E$ take a 
descending chain $(0,K_0),u_1,\ldots,u_n$ with $n\ge t^{1-\varepsilon}/2$ and thus 
$n\ge n_0$. Since $n\ge \lceil t^{1-\varepsilon}/2\rceil$ there exists $j\in\{1,\ldots,k\}$ 
and $l\in\{0,\ldots,n_0-m\}$ such that at least $m$ consecutive of the $n$ distances 
$d((0,K_0),u_1)$ and $d(u_{i-1},u_i)$, $2\le i\le n$, lie in the interval 
$I_j:=[s_{j-1},s_j]$. In addition
\begin{align*}
\|x_i\|\le\Big(\frac{t^{1-\varepsilon}}{2}+1\Big)\frac{2ct^\varepsilon}{4c}\le t,
\quad 1\le i\le n_0,
\end{align*}
where $x_i$ is the projection of $u_i$ onto the first coordinate, that is
$u_i=(x_i,K_i)$. Hence $E\subset \cup_{j=1}^k E_j$ where $\psi$ belongs to $E_j$
if it contains a descending chain $(x_1,K_1),\ldots,(x_m,K_m)$ with $\|x_1\|\le t$
and $d((x_{i-1},K_{i-1}),(x_i,K_i))\in I_j,\ i\in\{1,\ldots, m-1\}$. Hence
\begin{align}\label{ineq4.1}
\BP(\Psi^{(0,Z_0)}\in E)\le \sum_{j=1}^k \BP(\Psi\in E_j).
\end{align}
Furthermore we get
\begin{align*}
\BP(&\Psi\in E_j)\le \BE\Big[\sum_{((x_1,K_1),\ldots,(x_m,K_m))\in
\Psi^{(m)}} \I \{\|x_1\|\le t, d((x_i,K_i),(x_{i+1},K_{i+1}))\in I_j\}
\Big]\\
&=\int\cdots\int \I \{\|x_1\|\le t, d((x_i,K_i),(x_{i+1},K_{i+1}))\in
I_j\}
d x_1 \ldots d x_mQ(d K_1)\ldots Q(d K_m).
\end{align*}
A calculation similar to the one following \eqref{sum0} yields
\begin{align*}
\BP(\Psi\in E_j)&\le c_0 c_1^{m-1} t^d (s_j^d-s_{j-1}^d)^{m-1}
=c_0t^d\Big(\frac{c_1}{c_1+1}\Big)^{m-1}
\end{align*}
with a constant $c_0>0$. Hence there exist constants $c_2,c_3>0$ such that
\begin{align}\label{ineq4.2}
\BP(\Psi\in E_j)\le c_3 t^d \exp(-c_2 m),\quad j=1,\ldots,k.
\end{align}
Since $t^{\varepsilon d}\ge 1$ and $c>1$ we get $k\le (c_1+2) t^{\varepsilon d}$. 
Using \eqref{ineq4.1} and \eqref{ineq4.2} we obtain
\begin{align*}
\BP(\Psi^{(0,Z_0)}\in E)\le(c_2+2)t^{\varepsilon d} c_3 t^d \exp(-c_2 m)
\le c_4 t^{d+\varepsilon d} \exp\Big(-c_5
t^{1-\varepsilon-\varepsilon d}\Big)
\end{align*}
with constants $c_4,c_5>0$. Note that $1-\varepsilon-\varepsilon d >0$.
Now we give a bound for the first summand in \eqref{ineq4.3}.
By definition \eqref{defD} of $D=D(\Psi,0,Z_0)$ we have
\begin{align}\label{657}
\BP\left(D>\frac{t^\varepsilon}{4c}\right)
\le \BP\Big(\Psi\Big(B\Big(0,\frac{t^\varepsilon \zeta(Z_0)}{4c}\Big)
\times\mathcal K'_0\Big)=0\Big)
\le \BP\Big(\Psi\Big(B\Big(0,\frac{t^\varepsilon}{4c}\Big)\times\mathcal K'_0\Big)=0\Big),
\end{align}
where $\zeta(K), K\in\cK$, is the largest $r>0$ such that $B(0,r)\subset K$ .
The second inequality is due to the first inclusion in \eqref{shapecond}.
Since $\Psi$ is a Poisson process this yields
\begin{align*}
\BP\left(D>\frac{t^\varepsilon}{4c}\right)\le\exp(-c_6 t^{\varepsilon d}),
\end{align*}
where $c_6>0$.
With \eqref{ineq4.3} we get from the preceding two inequalities
\begin{align*}
\BP(U>t)\le \exp(-c_6 t^{\varepsilon d}) + c_4 t^{d+\varepsilon d} \exp\Big(-c_5
t^{1-\varepsilon-\varepsilon d}\Big),
\end{align*}
for all sufficiently large $t$. This implies the first assertion.

The second assertion can be shown using the same argumentation as above. 
But to be able to deal with the inserted point $(x,Z_0)$ we have to divide the 
event $E_j$ for $j=1,\ldots,k$ into events $E_{j,i}, i=1,\ldots,m$, where $E_{j,i}$ is 
the event that $E_j$ occurs and the $i$-th point of the descending chain of the 
configuration is $(x,K)$ for a fixed $K\in\cK_0$. Hence the analogue of \eqref{ineq4.1} 
consists of a double sum on the right-hand side. Since the arguments used to bound the 
probability of $E_j$ also hold for $E_{j,i}$ and because of the definition of $m$ 
we still obtain an sub-exponential decaying tail. We omit the details.
\qed

\bigskip

We define
\begin{align*}
&\rho_n(\psi,x,K,A):=\rho(n^{1/d}\psi,n^{1/d}x,K,n^{1/d}A)\I_{W_1}(x),\\
&\tau(s):=\sup_{n\ge 1,\ x\in W_1}
\BP(T(n^{1/d}((n^{-1/d}\Psi)_{W_1}-x),x,Z_0)>s),\quad s>0,\\
&\mu^\rho_n(A):=\sum_{(x,K)\in(n^{-1/d}\Psi)_{W_1}}\rho_n((n^{-1/d}\psi)_{W_1},x,K,A)
\end{align*}
for $\psi\in\bN_0$, $x\in\R^d$, $K\in\cK_0$ and $A\subset\R^d$. Note that
\begin{align}\label{distribeq}
\int_{W_n} f(n^{-1/d}x)\mu(\Psi_{W_n},dx)\stackrel{d}{=}\int_{W_1} f \ d\mu_n^\rho.
\end{align}
This equation links \eqref{cltstatement} with the form in which the central 
limit theorem is formulated in Theorem 2.2 in \cite{PeEJP}.
The next lemma shows that (2.7), (2.8) and (2.10) in \cite{PeEJP} are satisfied.

\begin{lemma}\label{the4.1}
Suppose $Z_1$ is an independent copy of $Z_0$ and also independent of $\Psi$. Then
\begin{align}\label{eqn4.1}
\BP(U(\Psi,0,Z_0)<\infty)=\BP(U(\Psi^{(x,Z_1)},0,Z_0)<\infty)=1.
\end{align}
Furthermore for $p>0$
\begin{align}
&\sup_{n\ge 1,\ x\in W_1}
\BE\Big[\rho_n\big((n^{-1/d}\Psi)_{W_1},x,Z_0,\R^d\big)^p\Big]<\infty,\label{eqn4.2}\\
&\sup_{n\ge 1,\ x, y\in W_1}
\BE\Big[\rho_n\big([(n^{-1/d}\Psi)_{W_1}]^{(y,Z_1)},x,Z_0,\R^d\big)^p\Big]<\infty,\label{eqn4.3}
\end{align}
and $\sup_{s\ge 1} s^q\tau(s)<\infty$ for some $q>d(150+6/p)$.
\end{lemma}
{\sc Proof:} The first assertion \eqref{eqn4.1} follows directly from Lemma \ref{l4.2}. 
For $p>0$ we have
\begin{align*}
\sup_{n\ge 1,\ y\in W_1}
\BE[\rho_n((n^{-1/d}\Psi)_{W_1},y,Z_0,\R^d)^p]&=\sup_{n\ge 1,\ y\in
W_1} \BE[\rho(n^{1/d}(n^{-1/d}\Psi)_{W_1},n^{1/d}y,Z_0,\R^d)^p]\\
&=\sup_{n\ge 1,\ x\in W_n} \BE[\rho(\Psi_{W_n},x,Z_0,\R^d)^p].
\end{align*}
The growth bound \eqref{gb} and  the definition of $\rho$ 
imply for all $n\ge 1$ and $x\in W_n$ that
\begin{align*}
\BE\big[\rho(\Psi_{W_n},x,Z_0,W_n)^p\big]\le
\alpha^p \BE\Big[R((\Psi_{W_n})^{(x,Z_0)},x,Z_0)^{\beta
p}\I_{\{{\rm{card}}((\Psi_{W_n})^{(x,Z_0)})>1\}}\Big].
\end{align*}
Define the random variable $Y_{n,x}:=
R((\Psi_{W_n})^{(x,Z_0)},x,Z_0)\I_{\{{\rm{card}}(\Psi_{W_n})>0\}}$. For $t\le\mbox{diam}(W_n)$ 
we conclude from $(2.10)$ in \cite{LaPe12}
\begin{align*}
\BP(Y_{n,x}>t)&\le \BP(R((\Psi_{W_n})^{(x,Z_0)},x,Z_0)>t)\le\exp(-c_0 t^d),
\end{align*}
where $c_0$ only depends on $W_1$. Trivially this inequality holds for $t>\mbox{diam}(W_n)$ 
as well. This yields \eqref{eqn4.2}. Assertion \eqref{eqn4.3} can be shown analogously, 
if we replace $Y_{n,x}$ by
\begin{align*}
Y_{n,x,y}:=R((\Psi_{W_n})^{(x,Z_0),(y,Z_1)},x,Z_0).
\end{align*}
To prove the statement about $\tau$ we note
\begin{align*}
\tau(s) =&\sup_{n\ge 1,\ y\in W_1} \BP(T(n^{1/d}(n^{-1/d}\Psi)_{W_1}-n^{1/d}y,Z_0)>s)\\
=&\sup_{n\ge 1,\ x\in W_n} \BP(T(\Psi_{W_n}-x,Z_0)>s)\\
\le &\sup_{n\ge 1,\ x\in W_n} \BP(U(\Psi_{W_n}-x,0,Z_0)>s),
\end{align*}
where we use $T(\varphi,K)\le U(\varphi,0,K)$ for $\varphi\in\bN_0,K\in\cK_0$.

With some minor technical modifications we can use the proof of Lemma \ref{l4.2} 
to show that $\BP(U(\Psi_{W_n}-x,0,Z_0)>s)$ is sub-exponentially decaying 
independent of $n\ge 1$ and $x\in W_n$. The dependence on $n$ vanishes because it 
is possible to refer to $\Psi$ instead of $\Psi_{W_n}$ for all $n\in\N$ if necessary. 
The dependence on $x$ vanishes because of the translation invariance of the 
Lebesgue measure. We omit the details.\qed

\bigskip

{\sc Proof of Theorem \ref{clt}:}
Because of Lemma \ref{the4.1} the assumptions of Theorem 2.1 and 2.2 in
\cite{PeEJP} are satisfied. Due to \eqref{distribeq} Theorem \ref{clt}
is implied by these two theorems.\qed

\bigskip

\begin{remark}\rm The first inclusion in the assumption \eqref{shapecond}
has been used to achieve the bound in \eqref{657}.
As the proof shows, this bound
can be replaced by the slightly weaker assumption
\begin{align*}
\BP(tB^d\subset Z_0)\ge 1-\exp(-(1/t)^\varepsilon)
\end{align*}
for all sufficiently small $t>0$.
\end{remark}

Consider a measurable function $h:\R_+\times\cK_0\rightarrow\R_+$.
If we choose $g(t,x,K,A)=h(t,K)\delta_x(A)$ for $t\ge0$, $x\in\R^d$,
$K\in\cK_0$, and measurable $A\subset\R^d$ we get
\begin{align*}
\mu(\psi,A)=
\I\{\card(\psi_{W_n})\ge 2\}\sum_{(x,K)\in\psi} h(R(\psi,x,K),K)\delta_x(A)),\quad \psi\in\bN_0.
\end{align*}
Recall that $V_d(\cdot)$ denotes the volume and assume, for instance,
that $h(t,K)=V_d(tK)$ and a constant density $f\equiv 1$. Then we can use
the hard-core property to deduce
\begin{align*}
\int_{W_n}f(n^{-1/d}x)\mu(\Psi_{W_n},dx)=V_d(Z(\Psi_{W_n}))\I\{\mbox{card}(\Psi_{W_n})\ge 2\},
\end{align*}
where $Z(\Psi_{W_n})$ is the union of the grains of the growth-maximal hard-core
model based on $\Psi_{W_n}$ as defined at \eqref{Z}. Similar considerations
lead to a central limit statement for the surface area in the case of strictly
convex particles.

\bigskip
\begin{center}\textbf{Acknowledgements}\end{center}
The authors thank Prof. Daryl Daley for numerous discussions on different
aspects of this work.
Moreover the authors are grateful to the DFG (Deutsche Forschungsgemeinschaft)
for partially funding this work via its research unit
``Geometry and Physics of Random Spatial Systems''.

\end{document}